\RequirePackage{fix-cm}
\documentclass[smallextended]{svjour3}       
\smartqed  
\usepackage{xspace}
\usepackage{algorithm}
\usepackage{subfigure}
\usepackage{amsfonts}
\usepackage{amsmath}
\usepackage{enumerate}
\usepackage{framed}
\usepackage{graphicx}
\usepackage{lscape}
\usepackage{bm}
\usepackage{color}
\usepackage{multicol}
\usepackage{cite}

\newcommand{\figref}[1]{{Figure~\ref{#1}}}

\newcommand{\thmref}[1]{{Theorem~\ref{#1}}}
\newcommand{\coref}[1]{{Corollary~\ref{#1}}}
\newcommand{\propref}[1]{{Proposition~\ref{#1}}}
\renewcommand{\eqref}[1]{{(\ref{#1})}}
\newcommand{\assref}[1]{{Assumption~\ref{#1}}}

 \newcommand{\N}{\mathbb{N}}

\newtheorem{dfn}{Definition}
\newtheorem{thm}{Theorem}
\newtheorem{lem}{Lemma}
\newcommand{\lemref}[1]{{Lemma~\ref{#1}}}
\newtheorem{ass}{Assumption}
\newtheorem{prp}{Proposition}
\newtheorem{rem}{Remark}
\newtheorem{cor}{Corollary}
\usepackage{graphicx}
\begin{document}

\title{Optimal strong convergence rates of some Euler-type timestepping schemes for the finite element discretization SPDEs driven by additive fractional Brownian motion and Poisson random measure}


\titlerunning{Optimal strong convergence of some schemes for SPDEs driven by additive fBm}        

\author{Aurelien Junior Noupelah, Antoine Tambue.
}

\authorrunning{A. J.  Noupelah, A. Tambue} 

\institute{A. J.  Noupelah \at
            Department of Mathematics and Computer Sciences, University of Dschang, P.O. BOX 67, Dschang, Cameroon\\
          \email{noupsjunior@yahoo.fr}   \\
          \and
          A. Tambue (Corresponding author) \at
          Department of Computing Mathematics and Physics,  Western Norway University of Applied Sciences, Inndalsveien 28, 5063 Bergen.\\           
            Center for Research in Computational and Applied Mechanics (CERECAM), and Department of Mathematics and Applied Mathematics, University of Cape Town, 7701 Rondebosch, South Africa.\\
             The African Institute for Mathematical Sciences(AIMS) of South Africa and Stellenbosh University,\\
             Tel.: +27-785580321\\
   \email{antonio@aims.ac.za, tambuea@gmail.com} 
}
\date{Received: date / Accepted: date}

\maketitle

\begin{abstract}
In this paper, we study the numerical approximation of a general second order semilinear stochastic partial differential equation (SPDE) driven by a additive 
fractional Brownian motion (fBm) with Hurst parameter $H>\frac 12$  and Poisson random measure, more realistic in modelling real world phenomena. 
To the best of our knowledge,  numerical schemes  for such SPDE have been lacked in scientific literature.
The approximation is done with the  standard finite element method  in space and  three Euler-type timestepping  methods in time, more precisely
 linear implicit method, exponential integrator and exponential Rosenbrock scheme are used for time discretisation.   
 In contract to the current literature in the field for SPDE driven only  by fBm,  our  linear operator is  not necessary self-adjoint and  optimal  strong convergence rates have been achieved for  SPDE driven  only by fBm and  SPDE driven by fBm and  Poisson measure.
  The results examine how the convergence orders depend on the regularity of the
 noise and the initial data and reveal that  the full discretization attains an optimal convergence rate of order $\mathcal{O}(h^2+\Delta t)$  for the exponential integrator and implicit schemes (linear operator $A$ self-adjoint for implicit). 
 Numerical experiments are provided to illustrate our theoretical results for the case of SPDE driven with fBm noise.

  
\end{abstract}

\section{Introduction}
\label{intro}
We analyse the strong numerical approximation of an SPDE defines in $\Lambda\subset\mathbb{R}^d$, $d \in \{1,2,3\}$ with initial value and boundary conditions (Dirichlet, Neumann, Robin boundary conditions or mixed Dirichlet and Neumann). 
In Hilbert space, our model equation can be formulated as the following   parabolic SPDE 
\begin{eqnarray}
 \label{model3}
 dX(t)=[AX(t)+F(X(t))]dt+\phi(t) dB^H(t)+\int_{\chi}z_0\widetilde{N}(dz,dt),\ X(0)=X_0
  \end{eqnarray}
  in Hilbert space $\mathcal{H}=L^2(\Lambda)$, with $\quad z_0\in\chi$, where  $\chi$ is the mark set defined by $\chi:=\mathcal{H}\setminus\{0\}$.   Let $\mathcal{B}(\Gamma)$ be the smallest $\sigma$-algebra containing all open sets of $\Gamma$. Let $(\chi, \mathcal{B}(\chi), \nu)$ be a $\sigma$-finite measurable space and $\nu$ ( with $\nu\not\equiv 0$) a L\'{e}vy measure on $\mathcal{B}(\chi)$ such that
\begin{eqnarray}
\nu(\{0\})=0\quad \text{and}\quad \int_{\chi}\min(\Vert z\Vert^2, 1)\nu(dz)<\infty.
\end{eqnarray}
 Let $N(dz, dt)$ be the $\mathcal{H}$-valued  Poisson  distributed $\sigma$-finite measure on the product $\sigma$-algebra $\mathcal{B}(\chi)$ and $\mathcal{B}(\mathbb{R}_{+})$ with intensity $\nu(dz)dt$, where $dt$ is the Lebesgue measure on $\mathcal{B}(\mathbb{R}_{+})$. In our model problem \eqref{model3}, $\widetilde{N}(dz, dt)$ stands for the compensated Poisson random measure defined by
 \begin{eqnarray}
 \widetilde{N}(dz, dt):=N(dz, dt)-\nu(dz)dt.
\end{eqnarray} 
 We denote by $T>0$, the final time, $F:\mathcal{H}\rightarrow \mathcal{H}$, $\phi$ are deterministic mappings that will be specified precisely later, $X_0$ is the initial data which is random, $-A$ is a linear operator, not necessary self-adjoint, unbounded  and generator of an analytic semigroup $S(t):=e^{-tA}$, $t\geq 0$. Note that $B^H(t)$ is a $\mathcal{H}$-valued $Q$-cylindrical fractional Brownian motion of Hurst parameter $H\in(\frac 12,1]$ in a filtered probability space $(\Omega, \mathcal{F}, \mathbb{P},\{\mathcal{F}_t\}_{t\geq 0})$ with the covariance operator $Q: \mathcal{H}\rightarrow \mathcal{H}$, which is  positive definite and  self-adjoint. 
The filtered probability space $(\Omega, \mathcal{F}, \mathbb{P},\{\mathcal{F}_t\}_{t\geq 0})$ is assumed to fulfil the usual condition 
(see \cite[Def 2.2.11]{Pre}). It is well known \cite{Car} that the noise can be represented as
\begin{eqnarray}
\label{def_fBm}
B^H(t)=\sum_{i\in\N^d}\beta^H_i(t)Q^{\frac 12}e_i=\sum_{i\in\N^d}\sqrt{q_i}\beta^H_i(t)e_i,\,\,\,t\in[0,T].
\end{eqnarray}
where $q_i,\,e_i,\,i\in\N^d$ are respectively the eigenvalues and eigenfunctions of the covariance operator $Q$, and $\beta_i^H$ are mutually independent and identically distributed fractional Brownian motions (fBm).

In  our study , we first study  in details  the following  particular case where $z_0=0$, i.e, the SPDE is driven only by fBm
\begin{eqnarray}
\label{pb1}
\left\{
  \begin{array}{ll}
    dX(t)+AX(t)dt=F(X(t))dt+\phi(t) dB^H(t), \\
    X(0)=X_0,\,\,\,\,t\in[0,T].
  \end{array}
\right.
\end{eqnarray} 
The self-similar and long-range dependence properties of the fBm make this process a suitable candidate to model  many phenomena  like financial markets (see, e.g., \cite{Che,Hu, Man}) and traffic networks (see, e.g., \cite{Lel,Wil}). 
In most cases, SPDEs of type \eqref{pb1} do not have explicit solutions and therefore  numerical algorithms are required for their approximations. It is important to  mention that  if $H\neq \frac 12$ the process $B^H$ is not a semi-martingale  and the standard stochastic calculus techniques are therefore obsolete while studying SPDEs of type \eqref{pb1}.  Alternative approaches to the standard It\^o calculus are therefore required in order to build a stochastic calculus framework for such fBm.  
In recent years, there have been various developments of stochastic calculus and stochastic differential equations with respect to the fBm especially for $H\in(\frac 12, 1]$ 
(see, for example \cite{Alo,Car,Mis})
and  theory of SPDEs driven by fractional Brownian motion  has been also studied. For example, linear and semilinear stochastic equations in a Hilbert space with an infinite dimensional fractional Brownian motion are considered in \cite{Duna,Dunb}. In contrast to standard Brownian ($H=1/2$) where there are numerous literature on numerical algorithms for SPDEs, few works have been done  for numerical methods for fBm for SPDEs of type  \eqref{pb1}.
Indeed, standard explicit and linear implicit schemes have been investigated in the literature for SPDEs of type \eqref{pb1} (see \cite{Kam,Kim,Wanc}). 
The works in \cite{Kam,Wanc} deal with self-adjoint operator and use the spectral Galerkin method for the spatial discretization.  This is very restrictive as many  concrete applications use  non self-adjoint operators. 
Beside numerical algorithms used for spatial discretization  and  time discretization in \cite{Kam,Wanc} are limited to few applications.
Our goal in this work  is to extend keys time stepping methods, which have been built for standard Brownian motion ($H=1/2$). These extensions are extremely complicated   due to  the fact that  the process $B^H$ is not a semi-martingale. 
Our results will be based on many novel intermediate lemmas. Indeed our  schemes here are based on finite element method (or finite volume method) for spatial discretization so that we gain the flexibility of 
these methods to deal with complex boundary conditions and we can apply well-developed techniques such as upwinding to deal with advection.
 For time discretization,  we first updated implicit linear for finite element method and not necessarily self-adjoint.  
 We also  provide the strong convergence of the exponential scheme \cite{Lor}  for ($H\in(\frac 12, 1]$).
  Note that this scheme is an explicit stable scheme, where the implementation is based on the computation of  matrix exponential  functions \cite{Lor}.
 As the linear implicit and exponential scheme are stable only when the linear operator  $A$ is stronger than the nonlinear function $F$ 
 \footnote{In this case the SPDE \eqref{pb1} is said to be driven by  its linear part.},  
 we also provide the  strong convergence  of the Stochastic Exponential Rosenbrock Scheme (SERS) \cite{Muka} for  ($H \in(\frac 12, 1]$), 
 which is very stable  when  \eqref{pb1} is driven both  by its linear or nonlinear part.
 
 However the model  equation  \eqref{pb1} can be unsatisfactory and  less realistic.  For instance, in finance, the unpredictable nature of many events 
 such as market crashes, announcements made by the central banks,  changing    credit risk, insurance in a changing risk, changing face of operational risk \cite{Tankov, Platen1} might have sudden and significant impacts on the stock price. 
  As for standard  Brownian motion,  we can  incorporate a non-Gaussian noise such as L\'{e}vy process or
   Poisson random measure to model such events.  
   The corresponding equation is our model equation given in \eqref{model3}. In contrast to SPDE driven by fBm in \eqref{pb1} where at least few numerical schemes exist, 
     numerical schemes  for such SPDE of type \eqref{model3} driven by fBm and Poisson measure  have been lacked in scientific literature, to the best of our knowledge.
     In  this work, we will also fill the gap by extending the implicit scheme, the exponential scheme and the Stochastic Exponential Rosenbrock Scheme to SPDE of type \eqref{model3}.
      For  SPDE  of type \eqref{pb1} and  SPDE of type \eqref{model3}, our strong convergence  results examine how the convergence orders depend on the regularity of the
 noise and the initial data and reveal that  the full discretization attains an optimal convergence rate of order $\mathcal{O}(h^2+\Delta t)$  for the exponential integrator and implicit schemes (linear operator $A$ self-adjoint for implicit).
%

The paper is structured as follows. In Section \ref{mathsetting}, Mathematical setting for fBm is presented, along with the well posedness and regularities results of the mild solution of SPDE \eqref{pb1} driven by fBm.  In Section \ref{schemes}, numerical schemes based on implicit scheme,  stochastic exponential integrator and stochastic  exponential Rosenbrock scheme for SPDE \eqref{pb1} driven by fBm are presented.   In Section \ref{conv}, the strong convergence proofs of schemes presented in Section \ref{schemes} are provided.   In Section \ref{convp},  numerical schemes based on semi implicit scheme,  stochastic exponential integrator  scheme and stochastic  exponential Rosenbrock scheme are presented for  SPDE \eqref{model3} driven by fBm and Poisson  measure, along with the extension of their strong convergence proofs. We end the paper in Section \ref{numerik} with numerical experiments illustrating our theoretical results  for SPDE \eqref{pb1} driven by fBm noise.

 \section{Mathematical setting}
 \label{mathsetting}
In this section, we review some standard  results on fractional calculus and introduce notations, definitions and preliminaries results which will be needed throughout this paper.
\begin{dfn}[\cite{Kam,Mas,Mis,Wanc}]
\label{deffBm}The fractional Brownian motion (fBm) of Hurst parameter $H\in (0,1)$ is a centered Gaussian  process $\beta^H=\{\beta^H(t),\,\,t\geq 0\}$ with the covariance function
\begin{eqnarray*}
\label{cov_fBm}
\mathbb{E}\left[\beta^H(t)\beta^H(s)\right]=\frac 12\left[t^{2H}+s^{2H}-|t-s|^{2H}\right].
\end{eqnarray*}
\end{dfn} 
Notice that if $H=\frac 12$, the process is the standard Brownian motion. 
\begin{rem}\cite[Remark 1.2.3]{Mis}
\label{deffBmH=1}
 For $H=1$, we set $\beta^H(t)=\beta^1(t)=t\xi$, where $\xi$ is a standard normal random variable.
\end{rem}
Throughout this paper  the Hurst parameter $H$  is assumed to be in the interval $(1/2,1]$.
Let $\left(K,\langle .,.\rangle_K,\|.\|\right)$ be a separable Hilbert space. For $p\geq 2$ and for a Banach space U, we denote by $L^p(\Omega,U)$ the Banach space of $p$-integrable $U$-valued random variables. 
We denote by $L(U,K)$ the space of bounded linear mapping from $U$ to $K$ endowed with the usual operator norm $\|.\|_{L(U,K)}$ and $\mathcal{L}_2(U,K)=HS(U,K)$ the space of Hilbert-Schmidt operators from $U$ to $K$ with 
\begin{eqnarray}
\label{norm1}
\left\|l\right\|_{\mathcal{L}_2(U,K)}:=\left(\sum_{i\in\N^d}\|le_i\|^2\right)^{\frac 12},\,\,\,\,\,l\in \mathcal{L}_2(U,K),
\end{eqnarray}
where $(e_i)_{i\in\N^d}$ is an orthonormal basis on $U$. The sum in $\eqref{norm1}$ is independent of the choice of the orthonormal basis in $U$. For simplicity, we use the notation $L(U,U)=:L(U)$ and $\mathcal{L}_2(U,U)=:\mathcal{L}_2(U)$. It is well known that for all $l\in L(U,K)$ and $l_1\in \mathcal{L}_2(U)$, $ll_1\in \mathcal{L}_2(U,K)$ and 
\begin{eqnarray}
\label{prop_norm1}
\|ll_1\|_{\mathcal{L}_2(U,K)}\leq \left\|l\right\|_{L(U,K)}\|l_1\|_{\mathcal{L}_2(U)}.
\end{eqnarray} 
We denote by $L^0_2$, the space of Hilbert-Schmidt operators from $Q^{\frac 12}(\mathcal{H})$ to $\mathcal{H}$ by $L^0_2:=HS(Q^{\frac 12}(\mathcal{H}),\mathcal{H})$ with corresponding norm $\|.\|_{L^0_2}$ by 
\begin{eqnarray}
\label{norm2}
\|l\|_{L_2^0}:=\left\|lQ^{\frac 12}\right\|_{HS}=\left(\sum_{i\in\N^d}\|lQ^{\frac 12}e_i\|^2\right)^{\frac 12},\,\,\,\,\,l\in L_2^0.
\end{eqnarray}
The following lemma will be very important throughout this paper.
\begin{lem}\cite[Lemma 1]{Car}
\label{lem1}
For any $\varphi\in L^{1/H}([0,T])$, the following inequality holds
\begin{eqnarray}
\label{maj}
H(2H-1)\int_{0}^{T}\int_{0}^{T}|\varphi(u)||\varphi(v)||u-v|^{2H-2}du dv\leq C_H\|\varphi\|^2_{L^{1/H}([0,T])}.
\end{eqnarray}
\end{lem}
Then for a function $f\in L^2\left([0,T],L^0_2\right)$, we consider the stochastic integral define as
\begin{eqnarray}
\label{stoint}
\int_0^Tf(s)dB^H(s):=\sum_{i\in\N^d}\int_0^Tf(s)Q^{\frac 12}e_id\beta_i^H(s).
\end{eqnarray} 
For all $i\in\mathbb{N}^d$, since the integrand $f$ is deterministic  then  the
 mean of the random variable $\int_0^Tf(s)d\beta_i^H(s)$ is zero and using 
 ( \cite[(2.12)]{Duna} or \cite[(2.14)]{Dunb}) 
 with \lemref{lem1}, the second moment (with $H\in(\frac 12,1)$) satisfies
\begin{eqnarray}
 \label{stoin}
 \mathbb{E}\left\|\int_0^Tf(s)d\beta_i^H(s)\right\|^2 \leq   C_H \|f\|^2_{L^{1/H}([0,T],\mathcal{H})},\hspace{1cm}H\in\left(1/2,1\right].
 \end{eqnarray}
As the sequence of random variables $\left(\int_0^Tf(s)Q^{\frac 12}e_id\beta_i^H(s),\,\,i\in\N^d\right)$ are mutually independent Gaussian random variable, the mean of random variable \eqref{stoint} is also zero and by \eqref{stoin}, we prove that its second moment satisfies
\begin{eqnarray}
\label{stoint1}
\mathbb{E}\left\|\int_0^Tf(s) dB^H(s)\right\|^2	&=&\mathbb{E}\left\|\sum_{i\in\N^d}\int_0^Tf(s)Q^{\frac 12}e_i d\beta_i^H(s)\right\|^2\nonumber\\
&=&\sum_{i\in\N^d}\mathbb{E}\left\|\int_0^Tf(s) Q^{\frac 12}e_i d\beta_i^H(s)\right\|^2\nonumber\\
&\leq& C_H\sum_{i\in\N^d}\left\|f(\cdot) Q^{\frac 12}e_i\right\|^2_{L^{\frac 1H}([0,T],\mathcal{H})}\nonumber\\
&=&C_H\sum_{i\in\N^d}\left(\int_0^T\left\|f(s) Q^{\frac 12}e_i\right\|^{\frac 1H}ds\right)^{2H},
\end{eqnarray}
and due to the fact that $\frac 1H<2$ then $L^2([0,T],\mathcal{H})$ is continuously embedded in $L^{\frac 1H}([0,T],\mathcal{H})$. Hence  for $f\in L^2\left([0,T],L^0_2\right)$
\begin{eqnarray}
\label{bound}
\mathbb{E}\left\|\int_0^Tf(s) dB^H(s)\right\|^2&\leq& C_H\sum_{i\in\N^d}\left\|f(\cdot) Q^{\frac 12}e_i\right\|^2_{L^{\frac 1H}([0,T],\mathcal{H})}\nonumber\\
&\leq& C_H\sum_{i\in\N^d}\left\|f(\cdot) Q^{\frac 12}e_i\right\|^2_{L^2([0,T],\mathcal{H})}\nonumber\\
&=&C_H\left\|f\right\|^2_{L^2([0,T],L^0_2)}<\infty.
\end{eqnarray}
In that follows, we will make some assumptions on $F$, $\phi$, $X_0$ and $A$, which will allow us to ensure the existence and uniqueness of the mild solution $X$ of \eqref{pb1} represented by (see e.g \cite{Wanc})
\begin{eqnarray}
\label{mild_sol_pb1}
X(t)=S(t)X_0+\int_{0}^{t}S(t-s)F(X(s))ds+\int_{0}^{t}S(t-s)\phi(s) dB^H(s),
\end{eqnarray}
for $t\in[0,T]$.
To ensure the existence and the uniqueness of solution for SPDE \eqref{pb1} and for the purpose of convergence analysis, we make the following assumptions.
\begin{ass}[Noise term]
\label{noise}
We assume that for some constant $\beta\in(0,1]$ and $\delta\in\left[\frac{2H+\beta-1}2,1\right]$, the deterministic mapping $\phi:[0,T]\times\Lambda\rightarrow L_2^0$ satisfies
\begin{eqnarray}
\label{ass_noise_term2}
\left\|A^{\frac{\beta-1}2}\phi(t)\right\|_{L_2^0}\leq C<\infty,\quad t\in[0,T]\\
\label{ass_noise_term3}
\left\|A^{\frac{\beta-1}2}(\phi(t_2)-\phi(t_1))\right\|_{L_2^0}\leq C_T(t_2-t_1)^{\delta},\quad 0\leq t_1\leq t_2\leq T.
\end{eqnarray}
\end{ass} 
\begin{ass}[Non linearity]
\label{non}
For the deterministic mapping $F:\mathcal{H}\rightarrow \mathcal{H}$, we assume that there exists constant $L\in(0,\infty)$ such that 
\begin{eqnarray}
\label{ass_non_lin1}
\|F(0)\|\leq L,\hspace{1cm}\|F(u)-F(v)\|\leq L\|u-v\|,\hspace{1cm}u,v\in \mathcal{H},
\end{eqnarray}
As a consequence of $\eqref{ass_non_lin1}$ it holds that
\begin{eqnarray}
\label{ass_non_lin3}
\|F(v)\|\leq L\left(1+\|v\|\right),\hspace{2cm}v\in \mathcal{H}.
\end{eqnarray}
\end{ass}
\begin{ass}[Initial Value]
\label{init}
We assume that $X_0:\Omega\rightarrow \mathcal{H}$ is a $\mathcal{F}_0/\mathcal{B}(\mathcal{H})$-measurable mapping and $X_0\in L^2\left(\Omega, D\left(A^{\frac{2H+\beta-1}{2}}\right)\right)$.
\end{ass}
In the Banach space $\mathcal{D}\left(A^{\frac{\alpha}{2}}\right)$, $\alpha\in\mathbb{R}$, we use notation $\left\|A^{\frac{\alpha}{2}}\cdot\right\|=\|\cdot\|_{\alpha}$ and 
we recall the following properties of the semigroup $S(t)$ generated by $-A$, that will be useful throughout this paper.
\begin{prp}[Smoothing properties of the semigroup]\cite{Pazy}
\label{semigroup}
Let $\alpha>0$, $\delta\geq 0$ and $0\leq \gamma\leq 1$, then there exist a constant $C>0$ such that 
\begin{eqnarray}
\label{semigroup_prp1}
\|A^{\delta}S(t)\|_{L(\mathcal{H})}\leq C t^{-\delta},\; \|A^{-\gamma}(I-S(t))\|_{L(\mathcal{H})}\leq C t^{\gamma},\quad t>0\\
\label{semigroup_prp2}
 \|D^l_tS(t)v\|\leq C t^{-l-(\gamma-\alpha)/2}\|v\|_{\alpha},\quad v\in D(A^{\alpha}).
\end{eqnarray}
where $l=0,1$ and $D^l=\frac{d^l}{dt^l}$. If $\delta>\gamma$ then $D(A^{\delta})\supset D(A^{\gamma})$. Moreover, $A^{\delta}S(t)=S(t)A^{\delta}\; \text{on}\; D(A^{\delta})$.
\end{prp}
The next lemma (specially \eqref{semigroup_prp5} and \eqref{semigroup_prp6}) is an important result which plays a crucial role to obtain regularity results, very useful in this work.
\begin{lem}
\label{lem2}
For any $0\leq\rho\leq 1$, $0\leq\gamma\leq 2$ and $0\leq \upsilon\leq H$ with $H\in \left(\frac 12,1\right]$, if the  linear operator is given by \eqref{operator}, there exists a positive constant $C$ such that for all $0\leq t_1\leq t_2\leq T$,
\begin{eqnarray}
\label{semigroup_prp3}
\int_{t_1}^{t_2}\|A^{\rho/2}S(t_2-r)\|^2_{L(\mathcal{H})}dr\leq C(t_2-t_1)^{1-\rho},\\
\label{semigroup_prp4}
\int_{t_1}^{t_2}\|A^{\gamma/2}S(t_2-r)\|_{L(\mathcal{H})}dr\leq C(t_2-t_1)^{1-\frac{\gamma}2},\\
\label{semigroup_prp5}
\int_{t_1}^{t_2}\|A^{H}S(t_2-r)\|^{\frac 1H}_{L(\mathcal{H})}dr\leq C,\\
\label{semigroup_prp6}
\int_{t_1}^{t_2}\|A^{\upsilon}S(t_2-r)\|^{\frac 1H}_{L(\mathcal{H})}dr\leq C(t_2-t_1)^{\frac{H-\upsilon}{H}}.
\end{eqnarray}
\end{lem}
\textit{Proof.} See \cite[Lemma 2.1]{Mukc} for the proof of \eqref{semigroup_prp3} and \eqref{semigroup_prp4}. Concerning the proof of \eqref{semigroup_prp5}, the border case $H=\frac 12$ if obtained using \eqref{semigroup_prp3} with $\rho=1$ and the order border case $H=1$ is also obtained using \eqref{semigroup_prp4} with $\gamma=2$. Hence the proof of \eqref{semigroup_prp5} is thus completed by interpolation theory. The proof of \eqref{semigroup_prp6} for $0\leq \upsilon\leq H$ is an immediate consequence
of Proposition \ref{semigroup}. The border case $\upsilon=H$ is proved by \eqref{semigroup_prp5}. This completes the proof of Lemma \ref{lem2}.$\hfill\square$
\begin{rem}
\label{disc_semigroup}
Proposition \ref{semigroup} and Lemma \ref{lem2} also hold with a uniform constant $C$ (independent of $h$) when $A$ and $S(t)$ are replaced respectively by their discrete versions $A_h$ and $S_h(t)$ defined in Section \ref{schemes}, see e.g. \cite{Lar,Lor}.
\end{rem}
The well posedness result is given in the following theorem along with optimal regularity
results in both space and time.
\begin{thm}
\label{reg}
Assume that Assumptions \ref{noise}-\ref{init} are satisfied, then there exists a unique mild solution given by $\eqref{mild_sol_pb1}$ such that for all $t\in[0,T]$, $X(t)\in L^2\left(\Omega,D\left(A^{\frac{2H+\beta-1}2}\right)\right)$ with 
\begin{eqnarray}
\label{mild_sol21}
\|X(t)\|_{L^2(\Omega,\mathcal{H})}\leq C\left(1+\|X_0\|_{L^2(\Omega,\mathcal{H})}\right),\\
 \|F(X(t))\|_{L^2(\Omega,\mathcal{H})}\leq C\left(1+\|X_0\|_{L^2(\Omega,\mathcal{H})}\right).
\end{eqnarray}
Moreover, if the linear operator is given by \eqref{operator}, the following optimal regularity results in space and time hold
\begin{eqnarray}
\label{mild_sol3}
\left\|A^{\frac{2H+\beta-1}2}X(t)\right\|_{L^2(\Omega,\mathcal{H})}\leq C \left(1+\left\|A^{\frac{2H+\beta-1}2}X_0\right\|_{L^2(\Omega,\mathcal{H})}\right),\quad t\in[0,T],
\end{eqnarray}
and for $0\leq t_1<t_2\leq T$;
\begin{eqnarray}
\label{mild_sol4}
\|X(t_2)-X(t_1)\|_{L^2(\Omega,\mathcal{H})}\leq C(t_2-t_1)^{\frac{2H+\beta-1}2} \left(1+\|A^{\frac{2H+\beta-1}2}X_0\|_{L^2(\Omega,\mathcal{H})}\right).
\end{eqnarray}
Where $C=C(\beta,L,T,H)$ is a positive constant and $\beta$ is the regularity parameter of Assumption \ref{noise}.
\end{thm}
\textit{Proof}
 \cite[Theorem 3.5]{Wanc} gives the result of existence and uniqueness of the mild solution $X$.
  For  regularity in space,  we adapt  from \cite[Theorem 2.1 (23), (24)]{Muka} by  just replacing $\beta$ in their case by $2H+\beta-1$. The difference will therefore be made at the level of the estimate of the stochastic integral
  \begin{eqnarray*}
  I^2=\mathbb{E}\left\|\int_0^tA^{\frac{2H+\beta-1}2}S(t-s)\phi(s) dB^H(s)\right\|^2.
  \end{eqnarray*}
To reach our goal,  we use triangle inequality, the estimate $(a+b)^2\leq 2a^2+2b^2$, \eqref{stoint1} and \eqref{bound}, \assref{noise}, \propref{semigroup}, \lemref{lem2} \eqref{semigroup_prp5} to have 
\begin{eqnarray}
\label{sto_integral}
I^2&=&\mathbb{E}\left\|\int_0^tA^{\frac{2H+\beta-1}2}S(t-s)\phi(s) dB^H(s)\right\|^2\nonumber\\
&\leq& 2\mathbb{E}\left\|\int_0^tA^{\frac{2H+\beta-1}2}S(t-s)\left(\phi(t)-\phi(s)\right) dB^H(s)\right\|^2\nonumber\\
&+&2\mathbb{E}\left\|\int_0^tA^{\frac{2H+\beta-1}2}S(t-s)\phi(t) dB^H(s)\right\|^2\nonumber\\
&\leq& 2C\int_0^t\left\|A^{\frac{2H+\beta-1}2}S(t-s)\left(\phi(t)-\phi(s)\right)\right\|^2_{L^0_2} ds\nonumber\\
&&+ 2C_H\sum_{i\in\N^d}\left(\int_0^t\left\|A^{\frac{2H+\beta-1}2}S(t-s)\phi(t) Q^{\frac 12}e_i\right\|^{\frac 1H}ds\right)^{2H}\nonumber\\
&\leq& 2C\int_0^t\left\|A^H S(t-s)\right\|^2_{L(\mathcal{H})}\left\|A^{\frac{\beta-1}2}\left(\phi(t)-\phi(s)\right)\right\|^2_{L^0_2} ds\nonumber\\
&&+ 2C_H\sum_{i\in\N^d}\left(\int_0^t\left\|A^{H}S(t-s)\right\|^{\frac 1H}_{L(\mathcal{H})}\left\|A^{\frac{\beta-1}2}\phi(t) Q^{\frac 12}e_i\right\|^{\frac 1H}ds\right)^{2H}\nonumber\\
&\leq& C\int_0^t(t-s)^{-2H}(t-s)^{2\delta}ds+ C_H\left(\sum_{i\in\N^d}\left\|A^{\frac{\beta-1}2}\phi(t) Q^{\frac 12}e_i\right\|^2\right)\nonumber\\
&&\times\left(\int_0^t\left\|A^{H}S(t-s)\right\|^{\frac 1H}_{L(\mathcal{H})}ds\right)^{2H}\nonumber\\
&\leq& C\int_0^t(t-s)^{2\delta-2H}ds+ C_H\left\|A^{\frac{\beta-1}2}\phi(t)\right\|^2_{L^0_2}\left(\int_0^t\left\|A^{H}S(t-s)\right\|^{\frac 1H}_{L(\mathcal{H})}ds\right)^{2H}\nonumber\\
&\leq& C t^{2\delta-2H+1}+ C\left\|A^{\frac{\beta-1}2}\phi(t)\right\|^2_{L^0_2}\leq C
\end{eqnarray}
For  the proof of \eqref{mild_sol4}, triangle inequality yields 
\begin{eqnarray*}
\|X(t_2)-X(t_1)\|_{L^2(\Omega,\mathcal{H})}&\leq& \left\|\left(S(t_2-t_1)-I\right)X(t_1)\right\|_{L^2(\Omega,\mathcal{H})}\nonumber\\
&+& \left\|\int_{t_1}^{t_2}S(t_2-s)F(X(s))ds\right\|_{L^2(\Omega,\mathcal{H})}\\
&+&\left\|\int_{t_1}^{t_2}S(t_2-s)\phi(s) dB^H(s)\right\|_{L^2(\Omega,\mathcal{H})}.
\end{eqnarray*}  
Using the stability property of the semigroup $S(t)$ \eqref{semigroup_prp1}  with $\gamma=\frac{2H+\beta-1}2$ and \eqref{mild_sol3} 
allows to have

\begin{eqnarray}
\label{proofreg}
&&\|X(t_2)-X(t_1)\|_{L^2(\Omega,\mathcal{H})}\nonumber\\
&\leq& \left\|A^{-\frac {2H+\beta-1}2}(S(t_2-t_1)-I)\right\|_{L(\mathcal{H})}\left\|A^{\frac{2H+\beta-1}2}X(t_1)\right\|_{L^2(\Omega,\mathcal{H})}\nonumber\\
&&+\int_{t_1}^{t_2}\left\|S_h(t_2-s)\right\|_{L(\mathcal{H})}\left\|F(X(s))\right\|_{L^2(\Omega,\mathcal{H})}ds\nonumber\\
&&+\left\|\int_{t_1}^{t_2}S(t_2-s)\phi(s) dB^H(s)\right\|_{L^2(\Omega,\mathcal{H})}\nonumber\\
&\leq& C(t_2-t_1)^{\frac{2H+\beta-1}2}\left(1+\|A^{\frac{2H+\beta-1}2}X_0\|_{L^2(\Omega,\mathcal{H})}\right)\nonumber\\
&&+C(t_2-t_1)\left(1+\|X_0\|_{L^2(\Omega,\mathcal{H})}\right)+\left\|\int_{t_1}^{t_2}S(t_2-s)\phi(s) dB^H(s)\right\|_{L^2(\Omega,\mathcal{H})}\nonumber\\
&\leq& C(t_2-t_1)^{\frac{2H+\beta-1}2}\left(1+\|A^{\frac{2H+\beta-1}2}X_0\|_{L^2(\Omega,\mathcal{H})}+\|X_0\|_{L^2(\Omega,\mathcal{H})}\right)\nonumber\\
&&+\left\|\int_{t_1}^{t_2}S(t_2-s)\phi(s) dB^H(s)\right\|_{L^2(\Omega,\mathcal{H})}\nonumber\\
&\leq& C(t_2-t_1)^{\frac{2H+\beta-1}2}\left(1+\|A^{\frac{2H+\beta-1}2}X_0\|_{L^2(\Omega,\mathcal{H})}\right)+II,
\end{eqnarray}
because $D\left(A^{\frac{2H+\beta-1}{2}}\right)$ is continuously embedded in $L^2\left(\Omega\right)$ and 
\begin{eqnarray*}
II:=\left\|\int_{t_1}^{t_2}S(t_2-s)\phi(s) dB^H(s)\right\|_{L^2(\Omega,\mathcal{H})}.
\end{eqnarray*}
For the estimate of $II$, using triangle inequality, the estimate $(a+b)^2\leq 2a^2+2b^2$, \eqref{stoint1} and \eqref{bound}, inserting an appropriate power of $A$, \propref{semigroup}, \assref{noise}, \lemref{lem2} \eqref{semigroup_prp6} with $\upsilon=\frac{1-\beta}2\in [0,\frac 12)$ (hence $0\leq \upsilon < H$), we obtain
\begin{eqnarray}
\label{stoint2}
II^2&=&\mathbb{E}\left\|\int_{t_1}^{t_2}S(t_2-s)\phi(s) dB^H(s)\right\|^2\nonumber\\
&\leq&2\mathbb{E}\left\|\int_{t_1}^{t_2}S(t_2-s)\left(\phi(t_2)-\phi(s)\right) dB^H(s)\right\|^2\nonumber\\
&+&2\mathbb{E}\left\|\int_{t_1}^{t_2}S(t_2-s)\phi(t_2) dB^H(s)\right\|^2\nonumber\\
&\leq& 2C\int_{t_1}^{t_2}\left\|S(t_2-s)\left(\phi(t_2)-\phi(s)\right)\right\|^2_{L^0_2} ds\nonumber\\
&+&2C_H\sum_{i\in\N^d}\left(\int_{t_1}^{t_2}\left\|S(t_2-s)\phi(t_2) Q^{\frac 12}e_i\right\|^{\frac 1H} ds\right)^{2H}\nonumber\\
&\leq& 2C\int_{t_1}^{t_2}\left\|S(t_2-s)A^{\frac{1-\beta}2}\right\|^2_{L(\mathcal{H})}\left\|A^{\frac{\beta-1}2}\left(\phi(t_2)-\phi(s)\right)\right\|^2_{L^0_2} ds\nonumber\\
&&+2C_H\sum_{i\in\N^d}\left(\int_{t_1}^{t_2}\left\|S(t_2-s)A^{\frac{1-\beta}2}\right\|^{\frac 1H}_{L(\mathcal{H})}\left\|A^{\frac{\beta-1}2}\phi(t_2) Q^{\frac 12}e_i\right\|^{\frac 1H} ds\right)^{2H}\nonumber\\
&\leq& C\int_{t_1}^{t_2}(t_2-s)^{\beta-1+2\delta}ds+2C_H\left(\sum_{i\in\N^d}\left\|A^{\frac{\beta-1}2}\phi(t_2) Q^{\frac 12}e_i\right\|^2\right)\nonumber\\
&&\times\left(\int_{t_1}^{t_2}\left\|S(t_2-s)A^{\frac{1-\beta}2}\right\|^{\frac 1H}_{L(\mathcal{H})} ds\right)^{2H}\nonumber\\
&\leq& C(t_2-t_1)^{\beta+2\delta}+C\left\|A^{\frac{\beta-1}2}\phi(t_2)\right\|^2_{L^0_2}\left(t_2-t_1\right)^{2H+\beta-1}\nonumber\\
&\leq& C(t_2-t_1)^{2H+2\beta-1}+C\left(t_2-t_1\right)^{2H+\beta-1}\leq C \left(t_2-t_1\right)^{2H+\beta-1}.
\end{eqnarray} 
Substituting \eqref{stoint2} in \eqref{proofreg} completes the proof of \eqref{mild_sol4} and therefore that of \thmref{reg}.$\hfill\square$
\section{Numerical schemes}
\label{schemes}
Throughout this section, we assume that $\Lambda$ is bounded and has smooth boundary or is a convex polygon of $\mathbb{R}^d$, $d\in \{1,2,3\}$. In the rest of this paper
we consider the  SPDE \eqref{pb1} to be of the following form
\begin{eqnarray}
\label{secondorder}
dX(t,x)+[-\nabla \cdot \left(\mathbf{D}\nabla X(t,x)\right)+\mathbf{q} \cdot \nabla X(t,x)]dt&=&f(x,X(t,x))dt\nonumber\\
&+& b(x,t) dB^H(t,x),
\end{eqnarray}
$ x\in\Lambda$, $t\in[0,T]$,
where the function $f : \Lambda\times \mathbb{R}\longrightarrow \mathbb{R}$ is continuously twice differentiable and the function
$b: \Lambda \times \mathbb{R}\longrightarrow \mathbb{R}$ is globally Lipschitz with respect to the second variable.
In the abstract framework \eqref{pb1}, the linear operator $A$ takes the form
\begin{eqnarray}
\label{operator}
Au=-\sum_{i,j=1}^{d}\dfrac{\partial}{\partial x_i}\left(D_{ij}(x)\dfrac{\partial u}{\partial x_j}\right)+\sum_{i=1}^dq_i(x)\dfrac{\partial u}{\partial x_i},\quad
\mathbf{D}=\left(D_{i,j} \right)_{1\leq i,j \leq d},\,\,
\end{eqnarray}
$\mathbf{q}=\left( q_i \right)_{1 \leq i \leq d}$, 
where $D_{ij}\in L^{\infty}(\Lambda)$, $q_i\in L^{\infty}(\Lambda)$. We assume that there exists a positive constant $c_1>0$ such that 
\begin{eqnarray}
\label{ellipticity}
\sum_{i,j=1}^dD_{ij}(x)\xi_i\xi_j\geq c_1|\xi|^2, \quad \forall \xi\in \mathbb{R}^d,\quad x\in\overline{\Omega}.
\end{eqnarray}
The functions $F : \mathcal{H} \longrightarrow \mathcal{H}$  and $ \phi : \mathbb{R}\longrightarrow HS\left(Q^{1/2}(\mathcal{H}), \mathcal{H}\right)$ are defined by 
\begin{eqnarray}
\label{nemystskii}
\left(F(v)\right)(x)=f\left(x,v(x)\right) 
\quad \text{and} \quad \left(\phi(t) (u)\right)(x)=b\left(x,t)\right).u(x),
\end{eqnarray}
for all $x\in \Lambda$, $v\in \mathcal{H}$, $u\in Q^{1/2}(\mathcal{H})$, with $\mathcal{H}=L^2(\Lambda)$.
 For an appropriate family of eigenfunctions $(e_i)$ such that $\sup\limits_{i\in\mathbb{N}^d}\left[\sup\limits_{x\in \Lambda}\Vert e_i(x)\Vert\right]<\infty$, 
 it is well known \cite[Section 4]{Arnulf1} that the Nemytskii operator $F$ related to $f$
 and the  operator $\phi$ associated  to $b$ defined in \eqref{nemystskii} satisfy Assumption \ref{noise} and Assumption \ref{non}.
As in \cite{Lor,Fuj} we introduce two spaces $\mathbb{H}$ and $V$, such that $\mathbb{H}\subset V$; the two spaces depend on the boundary 
conditions and the domain of the operator $A$. For  Dirichlet (or first-type) boundary conditions we take 
\begin{eqnarray*}
V=\mathbb{H}=H^1_0(\Lambda)=\overline{C^{\infty}_{c}(\Lambda)}^{H^1(\Lambda)}=\{v\in H^1(\Lambda) : v=0\quad \text{on}\quad \partial \Lambda\}.
\end{eqnarray*}
For Robin (third-type) boundary condition and  Neumann (second-type) boundary condition, which is a special case of Robin boundary condition, we take $V=H^1(\Lambda)$
\begin{eqnarray*}
\mathbb{H}=\{v\in H^2(\Lambda) : \partial v/\partial \mathtt{v}_{ A}+\alpha_0v=0,\quad \text{on}\quad \partial \Lambda\}, \quad \alpha_0\in\mathbb{R},
\end{eqnarray*}
where $\partial v/\partial \mathtt{v}_{ A}$ is the normal derivative of $v$ and $\mathtt{v}_{ A}$ is the exterior pointing normal $n=(n_i)$ to the boundary of $A$, given by
\begin{eqnarray*}
\partial v/\partial\mathtt{v}_{A}=\sum_{i,j=1}^dn_i(x)D_{ij}(x)\dfrac{\partial v}{\partial x_j},\,\,\qquad x \in \partial \Lambda.
\end{eqnarray*}
Using the Green's formula and the boundary conditions, the  corresponding bilinear form associated to $A$  is given by
\begin{eqnarray*}
a(u,v)=\int_{\Lambda}\left(\sum_{i,j=1}^dD_{ij}\dfrac{\partial u}{\partial x_i}\dfrac{\partial v}{\partial x_j}+\sum_{i=1}^dq_i\dfrac{\partial u}{\partial x_i}v\right)dx, \quad u,v\in V,
\end{eqnarray*}
for Dirichlet and Neumann boundary conditions, and  
\begin{eqnarray*}
a(u,v)=\int_{\Lambda}\left(\sum_{i,j=1}^dD_{ij}\dfrac{\partial u}{\partial x_i}\dfrac{\partial v}{\partial x_j}+\sum_{i=1}^dq_i\dfrac{\partial u}{\partial x_i}v\right)dx+\int_{\partial\Lambda}\alpha_0 uvdx, \quad u,v\in V,
\end{eqnarray*}
for Robin boundary conditions. Using the G\aa rding's inequality, it holds that there exist two constants $c_0$ and $\lambda_0$ such that
\begin{eqnarray}
\label{ellip1}
a(v,v)\geq \lambda_0\Vert v \Vert^2_{H^1(\Lambda)}-c_0\Vert v\Vert^2, \quad v\in V.
\end{eqnarray}
By adding and substracting $c_0Xdt$ in both sides of \eqref{pb1}, we have a new linear operator
 still denoted by $A$, and the corresponding  bilinear form is also still denoted by $a$. Therefore, the following coercivity property holds
\begin{eqnarray}
\label{ellip2}
a(v,v)\geq \lambda_0\Vert v\Vert^2_1,\quad v\in V.
\end{eqnarray}
Note that the expression of the nonlinear term $F$ has changed as we included the term $c_0X$ in a new nonlinear term that we still denote by  $F$. The coercivity property (\ref{ellip2}) implies that $-A$ is sectorial in $L^{2}(\Lambda)$, i.e.  there exist $C_{1},\, \theta \in (\frac{1}{2}\pi,\pi)$ such that
\begin{eqnarray}
\label{prop_op}
 \Vert (\lambda I +A )^{-1} \Vert_{L(L^{2}(\Lambda))} \leq \dfrac{C_{1}}{\vert \lambda \vert },\;\quad \quad
\lambda \in S_{\theta},
\end{eqnarray}
where $S_{\theta}=\left\lbrace  \lambda \in \mathbb{C} :  \lambda=\rho e^{i \phi},\; \rho>0,\;0\leq \vert \phi\vert \leq \theta \right\rbrace $ (see \cite{Henry}).
 Then  $-A$ is the infinitesimal generator of a bounded analytic semigroup $S(t)=e^{-t A}$  on $L^{2}(\Lambda)$  such that
\begin{eqnarray}
\label{semigroup1}
S(t)= e^{-t A}=\dfrac{1}{2 \pi i}\int_{\mathcal{C}} e^{ t\lambda}(\lambda I +A)^{-1}d \lambda,\;\;\;\;\;\;\;
\;t>0,
\end{eqnarray}
where $\mathcal{C}$  denotes a path that surrounds the spectrum of $-A $.
The coercivity  property \eqref{ellip2} also implies that $A$ is a positive operator and its fractional powers are well defined  for any $\alpha>0,$ by
\begin{equation}
\label{fractional}
 \left\{\begin{array}{rcl}
         A^{-\alpha} & =& \frac{1}{\Gamma(\alpha)}\displaystyle\int_0^\infty  t^{\alpha-1}{\rm e}^{-tA}dt,\\
         A^{\alpha} & = & (A^{-\alpha})^{-1},
        \end{array}\right.
\end{equation}
where $\Gamma(\alpha)$ is the Gamma function (see \cite{Henry}).
Under condition \eqref{ellipticity}, it is well known (see e.g. \cite{Fuj}) that the linear operator $-A$ given by \eqref{operator} generates an analytic semigroup $S(t)\equiv e^{-tA}$. 
Following \cite{Lor,Fuj}, we characterize the domain of the operator $A^{r/2}$ denoted by $\mathcal{D}(A^{r/2})$, $r\in\{1,2\}$ with the following equivalence of norms, useful in our convergence proofs
\begin{eqnarray*}
\Vert v\Vert_{H^1(\Lambda)}\equiv \Vert A^{r/2}v\Vert=:\Vert v\Vert_r,\quad \hspace{2cm} \forall v\in\mathcal{D}(A^{r/2}),\\
\mathcal{D}(A^{r/2})=\mathbb{H}\cap H^{r}(\Lambda), \quad \text{ (for Dirichlet boundary conditions)},\\
\mathcal{D}(A)=\mathbb{H}, \quad \mathcal{D}(A^{1/2})=H^1(\Lambda), \quad \text{(for Robin boundary conditions)}.
\end{eqnarray*}

We consider the discretization of the spatial domain by a finite element triangulation \cite{Mukb,Wanb}.
Let $\mathcal{T}_h$ be a set of disjoint intervals of $\Omega$ (for $d=1$), a triangulation of $\Omega$ (for $d=2$) or a set of tetrahedra (for $d=3$) 
with maximal length $h$ satisfying the usual regularity assumptions.\\
Let $V_h\subset \mathcal{H}$ denote the space of continuous functions that are piecewise linear over triangulation $\mathcal{T}_h$. To discretise in space, we introduce the projection $P_h$ from $L^2(\Omega)$ to $V_h$ define for $u\in L^2(\Omega)$ by
\begin{eqnarray}
\label{proj}
\left\langle P_hu,\mathcal{X}\right\rangle =\left\langle u,\mathcal{X}\right\rangle ,\hspace{2cm}\forall \mathcal{X}\in V_h.
\end{eqnarray} 
The discrete operator $A_h:V_h\rightarrow V_h$ is defined by
\begin{eqnarray}
\label{discrete_op}
\left\langle A_h\varphi,\mathcal{X}\right\rangle =-a\left\langle \varphi,\mathcal{X}\right\rangle ,\hspace{2cm}\varphi,\mathcal{X}\in V_h,
\end{eqnarray}
where $a$ is the corresponding bilinear form of A.

Like the operator $A$, the discrete operator $-A_h$ is also the generator of an analytic semigroup $S_h(t):=e^{-tA_h}$. The semi-discrete space version of problem \eqref{pb1} is to find $X^h(t)=X^h(\cdot,t)$ such that for $t\in[0,T]$
\begin{eqnarray}
\label{pb2}
\left\{
  \begin{array}{ll}
    dX^h(t)+A_hX^h(t)dt=P_hF(X(t))dt+P_h\phi(t) dB^H(t), \\
    X^h(0)=P_hX_0,\hspace{2cm}t\in[0,T].
  \end{array}
\right.
\end{eqnarray}
The mild solution of \eqref{pb2} can be represented as follows
\begin{eqnarray}
\label{mild_sol_pb2}
X^h(t)&=&S_h(t)X^h(0)+\int_{0}^{t}S_h(t-s)P_hF(X^h(s))ds\nonumber\\
&+&\int_{0}^{t}S_h(t-s)P_h\phi(s) dB^H(s),
\end{eqnarray} 
and we have the following regularity results.
\begin{lem}
\label{lem3}
Assume that Assumptions \ref{noise}-\ref{init} are satisfied, then the unique mild solution $X^h(t)$ given by \eqref{mild_sol_pb2} satisfied
\begin{eqnarray}
\label{reg1}
\left\|A^{\frac{2H+\beta-1}2}X^h(t)\right\|_{L^2(\Omega,\mathcal{H})}\leq C \left(1+\left\|A^{\frac{2H+\beta-1}2}X_0\right\|_{L^2(\Omega,\mathcal{H})}\right)\, t\in[0,T],
\end{eqnarray}
and for $0\leq t_1<t_2\leq T$;
{\small
\begin{eqnarray}
\label{reg2}
\|X^h(t_2)-X^h(t_1)\|_{L^2(\Omega,\mathcal{H})}\leq C(t_2-t_1)^{\frac{2H+\beta-1}2} \left(1+\|A^{\frac{2H+\beta-1}2}X_0\|_{L^2(\Omega,\mathcal{H})}\right).
\end{eqnarray}
}
\end{lem}
\textit{Proof} Since the operators $A_h$ and $S_h(t)$ satisfy the same properties as $A$ and $S(t)$ ( see Remark \ref{disc_semigroup}), then by using \cite[(83)]{Mukb} and  the boundedness of $P_h$ in the proof of \eqref{mild_sol3} and \eqref{mild_sol4}, we obtain the proof of  the expression \eqref{reg1} and \eqref{reg2}. The proof of Lemma \ref{lem3} is thus completed. $\hfill\square$

Now applying the linear implicit Euler method \cite{Kam,Wanb} to \eqref{pb2}  gives the following fully discrete scheme
\begin{eqnarray}
\label{impl}
\left\{
  \begin{array}{ll}
  	X^h_0=P_hX_0\\
    X^h_{m+1}=S_{h,\Delta t}X^h_m+\Delta tS_{h,\Delta t}P_hF(X^h_m)\Delta t+S_{h,\Delta t}P_h\phi(t_m) \Delta B^H_m. 
  \end{array}
\right.
\end{eqnarray} 
Furthermore applying the stochastic exponential integrator (\cite{Lor}, SETD1) and Rosenbrock scheme (\cite{Muka}, SERS) to \eqref{pb2} yields
\begin{eqnarray}
\label{SETD1}
\left\{
  \begin{array}{ll}
  	Y^h_0=P_hX_0\\
    Y^h_{m+1}=S_h(\Delta t)\left(Y^h_m+P_h\phi(t_m)\Delta B^H_m\right)+\Delta t \varphi_1(\Delta tA_h)P_hF(Y^h_m), 
  \end{array}
\right.
\end{eqnarray} 
and
\begin{eqnarray}
\label{SERS}
\left\{
  \begin{array}{lll}
  	Z^h_0=P_hX_0\\
    Z^h_{m+1}=e^{(-A_h+J^h_m)\Delta t}Z^h_m+(-A_h+J^h_m)\left(e^{(-A_h+J^h_m)\Delta t}-I\right)G^h_m(Z^h_m)\nonumber\\
    \hspace{2cm}+e^{(-A_h+J^h_m)\Delta t}P_h\phi(t_m)\Delta B^H_m. 
  \end{array}
\right.
\end{eqnarray}   
Where  $\varphi_1(\Delta tA_h)=(\Delta tA_h)^{-1}\left(e^{\Delta tA_h}-I\right)=\frac 1{\Delta t}\int_{0}^{\Delta t}e^{(\Delta t-s)A_h} ds$, $J^h_m$ is the Frechet derivative of $P_hF$ at $Z^h_m$ and $S_{h,\Delta t}:=(I+\Delta t A_h)^{-1}$. The term $G^h_m$ is the remainder at $Z^h_m$ and defines for all $\omega\in\Omega$ by
\begin{eqnarray}
\label{fre_der}
J^h_m:=(P_hF)'(Y^h_m(\omega))=P_hF'(Y^h_m(\omega)),\\
\label{remain}
G^h_m(\omega)(Z^h_m):=P_hF(Z^h_m)-J^h_m(\omega)Z^h_m,
\end{eqnarray}
and 
\begin{eqnarray*}
\Delta B^H_m:=B^H_{m+1}-B^H_m=\sum_{i\in\N^d}\sqrt{q_i}(\beta^H_i(t_{m+1})-\beta^H_i(t_m))e_i.
\end{eqnarray*}
Note  that  the exponential integrator scheme  \eqref{SETD1} is an explicit stable scheme  when  the  SPDE \eqref{pb1} is driven by its linear part as the linear implicit method, while the Stochastic Exponential Rosenbrock Scheme (SERS) \eqref{SERS} is very stable  when  \eqref{pb1} is driven by its linear or nonlinear part. When dealing with SERS, the strong convergence proof will make use of the following assumption, also used in
\cite{Muka,Mukc}.
\begin{ass}
\label{der}
For the deterministic mapping $F:\mathcal{H}\rightarrow \mathcal{H}$, we also assume that there exists constant $L\in(0,\infty)$ such that 
\begin{eqnarray}
\label{ass_non_lin2}
\|F'(u)v\|\leq L\|v\|,\hspace{2cm}u,v\in \mathcal{H}.
\end{eqnarray}
\end{ass}

\subsection{Main Result for   SPDE  driven by fBm}
\begin{thm}
\label{cstrongconvthm}
Let $X(t_m)$ be the mild solution of \eqref{pb1} at time $t_m=m\Delta t$, $\Delta t\geq 0$ represented by \eqref{mild_sol_pb1}. Let $\zeta^h_m$ be the numerical approximations through \eqref{impl} and \eqref{SERS}($\zeta^h_m=X^h_m$ for implicit scheme, $\zeta^h_m=Z^h_m$ for SERS). Under Assumptions \ref{noise}-\ref{init} and \ref{der} (essentially for SERS), $\beta\in(0,1]$, then the following holds
\begin{eqnarray}
\label{strongconvSETD1}
\left(\mathbb{E}\|X(t_m)-Y^h_m\|^2\right)^{\frac 12}\leq C\left(h^{2H+\beta-1}+\Delta t^{\frac{\beta+2H-1}2}\right),
\end{eqnarray}
and 
\begin{eqnarray}
\label{strongconvimplandSERS}
\left(\mathbb{E}\|X(t_m)-\zeta^h_m\|^2\right)\leq C \left(h^{2H+\beta-1}+\Delta t^{\frac{\beta+2H-1}2-\epsilon}\right),
\end{eqnarray}
where $\epsilon$ is a positive constant small enough.
\end{thm}
\section{Proofs of the main result for SPDE with fBm}
\label{conv}
We introduce the Riesz representation $R_h:V\rightarrow V_h$ defined by 
\begin{eqnarray}
\label{Riez_op1}
\langle AR_hv,\mathcal{X}\rangle=\langle Av,\mathcal{X}\rangle=a(v,\mathcal{X}),\hspace{1cm}v\in V,\hspace{1cm}\forall\mathcal{X}\in V_h
\end{eqnarray}
under the regularity assumptions on the triangulation and in view of the V-ellipticity, it is well known (\cite{Lar,Fuj}) that the following error bound holds:
\begin{eqnarray}
\label{Riez_op2}
\|R_hv-v\|+h\|R_hv-v\|_{H^1(\Omega)}\leq Ch^r\|v\|_{H^r(\Omega)},\quad v\in V\cap H^r(\Omega),
\end{eqnarray}
for $ r\in[1,2]$.
Let us consider the following deterministic linear problem:\\
Find $u\in V$ such that
\begin{eqnarray}
\label{pb4}
\frac{du}{dt}+Au=0,\hspace{1cm}u(0)=v,\hspace{0.5cm}t\in(0,T].
\end{eqnarray}
The corresponding semi-discrete problem in space consists to finding $u^h\in V_h$ such that 
\begin{eqnarray}
\label{pb5}
\frac{du^h}{dt}+Au^h=0,\hspace{1cm}u^h(0)=P_hv,\hspace{0.5cm}t\in(0,T].
\end{eqnarray}
Let us define the following operator
\begin{eqnarray}
\label{op}
G_h(t):=S(t)-S_h(t)P_h=e^{-At}-e^{-A_ht}P_h
\end{eqnarray}
Then we have the following lemma 
\begin{lem} 
\label{lem4}
The following estimates hold  for  the semi-discrete approximation of \eqref{pb2}. There exists a constant $C>0$ such that 
\begin{enumerate}
\item [(i)] For $v\in D(A^{\gamma/2})$
\begin{eqnarray}
\label{spa_err1}
\|u(t)-u^h(t)\|=\|G_h(t)v\|\leq Ch^rt^{-(r-\gamma)/2}\|v\|_{\gamma},\; r\in[0,2],\;\gamma\leq r
\end{eqnarray}
for any $t\in(0,T]$. 
\item[(ii)]For $v\in D(A^{\frac{\gamma-1}2})$
\begin{eqnarray}
\label{spa_err2}
\left(\int_0^t\|G_h(s)v\|^2ds\right)^{\frac 12}\leq Ch^{\gamma}\|v\|_{\gamma-1},\hspace{1cm}0\leq \gamma\leq 2,\hspace{1cm}t>0.
\end{eqnarray}
\item[(iii)]For $v\in D(A^{-\frac{\rho}2})$
\begin{eqnarray}
\label{spa_err3}
\left\|\int_0^tG_h(s)vds\right\|\leq Ch^{2-\rho}\|v\|_{-\rho},\hspace{1cm}0\leq \rho\leq 1,\hspace{1cm}t>0.
\end{eqnarray}
\item[(iv)]For $v\in D(A^{\frac{\delta-1}2})$
\begin{eqnarray}
\label{spa_err4}
\left(\int_0^t\|G_h(s)v\|^{\frac 1H}ds\right)^{2H}\leq Ch^{2(2H+\delta-1)}\|v\|^2_{\delta-1},\, 0\leq \delta\leq 1,\, t>0.
\end{eqnarray}
\end{enumerate}
\end{lem}
\textit{Proof} (see \cite[Lemma 3.1 and Lemma 3.2 (iv) and (v)]{Mukb}) for the proof of $(i)-(iii)$.  Let us prove  $(iv)$.
\begin{itemize}
\item For $H=\frac 12$, using \eqref{spa_err2} with $\gamma=\delta$, we obtain 
\begin{eqnarray}
\label{spa_err41}
\left(\int_0^t\|G_h(s)v\|^{\frac 1H}ds\right)^{2H}&=&\int_0^t\|G_h(s)v\|^2ds\nonumber\\
&\leq& Ch^{2\delta}\|v\|^2_{\delta-1}=Ch^{2(2H+\delta-1)}\|v\|^2_{\delta-1}.
\end{eqnarray}
\item For $H=1$, using \eqref{spa_err3} with $\rho=1-\delta$, we obtain
\begin{eqnarray}
\label{spa_err42}
\left(\int_0^t\|G_h(s)v\|^{\frac 1H}ds\right)^{2H}&=&
\left(\int_0^t\|G_h(s)v\|ds\right)^2\nonumber\\
&\leq& C\left(h^{2-(1-\delta)}\|v\|_{\delta-1}\right)^2\nonumber\\
&=&Ch^{2+2\delta}\|v\|^2_{\delta-1}=Ch^{2(2H+\delta-1)}\|v\|^2_{\delta-1}.
\end{eqnarray}
\end{itemize}
 Hence the proof of \eqref{spa_err4} is thus completed by interpolation theory.$\hfill\square$
\begin{lem}[Space error]
\label{serror}
Let Assumptions \ref{noise}-\ref{init} be fulfilled, then the following error estimate holds for the mild solution \eqref{mild_sol_pb1} and the discrete problem \eqref{mild_sol_pb2} holds
\begin{eqnarray}
\label{spa_err}
\|X(t)-X^h(t)\|_{L^2(\Omega,\mathcal{H})}\leq Ch^{2H+\beta-1}.
\end{eqnarray}
\end{lem}
\textit{Proof} Using triangle inequality, we have 
\begin{eqnarray*}
e(t)&:=&\|X(t)-X^h(t)\|_{L^2(\Omega,\mathcal{H})}\nonumber\\
&=&\Vert S(t)X_0+\int_{0}^{t}S(t-s)F(X(s))ds+\int_{0}^{t}S(t-s)\phi(s) dB^H(s)\nonumber\\
&&-S_h(t)P_hX_0-\int_{0}^{t}S_h(t-s)P_hF(X^h(s))ds\nonumber\\
&&-\int_{0}^{t}S_h(t-s)P_h\phi(s) dB^H(s)\Vert_{L^2(\Omega,\mathcal{H})}\nonumber\\
&\leq& \|G_h(t)X_0\|_{L^2(\Omega,\mathcal{H})}\nonumber\\
&+&\left\|\int_{0}^{t}S(t-s)F(X(s))ds-\int_{0}^{t}S_h(t-s)P_hF(X(s))ds\right\|_{L^2(\Omega,\mathcal{H})}\nonumber\\
&& +\left\|\int_{0}^{t}G_h(t-s)\phi(s) dB^H(s)\right\|_{L^2(\Omega,\mathcal{H})}\nonumber\\
&=:&e_0(t)+e_1(t)+e_2(t),
\end{eqnarray*}
with
\begin{eqnarray*}
&&e_0(t)=\|G_h(t)X_0\|_{L^2(\Omega,\mathcal{H})},\\
&&e_1(t)=\left\|\int_{0}^{t}S(t-s)F(X(s))ds-\int_{0}^{t}S_h(t-s)P_hF(X^h(s))ds\right\|_{L^2(\Omega,\mathcal{H})}\\
\text{and}
&& e_2(t)=\left\|\int_{0}^{t}G_h(t-s)\phi(s) dB^H(s)\right\|_{L^2(\Omega,\mathcal{H})}.
 \end{eqnarray*}
 Note that the deterministic error is already estimated, so we will  mostly concentrate our study on the stochastic error. 
  Indeed  Lemma \ref{lem4} with $r=\gamma=2H+\beta-1$ yields
 \begin{eqnarray}
 \label{space_err0}
e_0(t)=\|G_h(t)X_0\|_{L^2(\Omega,V)}\leq C h^{2H+\beta-1}\|A^{\frac{2H+\beta-1}2}X_0\|_{L^2(\Omega,\mathcal{H})}.
 \end{eqnarray}
 Using triangle inequality, the boundedness of $S_h(t-s)$ and $P_h$, Assumption \ref{non} (more precisely \eqref{ass_non_lin1}), we estimate the error $e_1(t)$ as follow
\begin{equation}
\label{spa_err_1}
\begin{split}
e_1(t)&=\left\|\int_{0}^{t}S(t-s)F(X(s))ds-\int_{0}^{t}S_h(t-s)P_hF(X(s))ds\right\|_{L^2(\Omega,\mathcal{H})}\\
&\leq \left\|\int_{0}^{t}\left(S(t-s)-S_h(t-s)P_h\right)F(X(s))ds\right\|_{L^2(\Omega,\mathcal{H})}\\
&\,\,\,\,\,\,\,+\left\|\int_{0}^{t}S_h(t-s)P_h\left(F(X(s))-F(X^h(s))\right)ds\right\|_{L^2(\Omega,\mathcal{H})}\\
&\leq \left\|\int_{0}^{t}\left(S(t-s)-S_h(t-s)P_h\right)\left(F(X(t))-F(X(s))\right)ds\right\|_{L^2(\Omega,\mathcal{H})}\\
&\,\,\,\,\,\,\,+\left\|\int_{0}^{t}\left(S(t-s)-S_h(t-s)P_h\right)F(X(t))ds\right\|_{L^2(\Omega,\mathcal{H})}+C\int_{0}^{t}e(s)ds.
\end{split}
\end{equation}
Applying Lemma \ref{lem4} $(i)$ with $r=2H+\beta-1$ and $\gamma=0$, Assumption \ref{non} \eqref{ass_non_lin1}, Theorem \ref{reg} (more precisely \eqref{mild_sol4}) to the first term and Lemma \ref{lem4} $(iii)$ with $\rho=0$, Theorem \ref{reg} (more precisely \eqref{mild_sol21}) and Lemma \ref{lem3} 
 to the second term yields
\begin{equation}
\label{space_err11}
\begin{split}
e_1(t)&\leq Ch^{2H+\beta-1}\int_{0}^{t}(t-s)^{-\frac{2H+\beta-1}2}\left\|X(t)-X(s)\right\|_{L^2(\Omega,\mathcal{H})}ds\\
&\,\,\,\,\,\,\,\,+C h^2\left\|F(X(t))\right\|_{L^2(\Omega,\mathcal{H})}+C\int_{0}^{t}e(s)ds\\
&\leq  Ch^{2H+\beta-1}+C\int_{0}^{t}e(s)ds.
\end{split}
\end{equation}  
For the estimation of $e_2(t)$, triangle inequality, the estimate $(a+b)^2\leq 2a^2+2b^2$, \eqref{stoint1} and \eqref{bound}, \lemref{lem4} ((i) with $r=2H+\beta-1$, $\gamma=\beta-1$, $(iv)$ with $\delta=\beta$) and Assumption \ref{noise} yields

\begin{eqnarray}
\label{spa_err_2}
e_2(t)^2&=& \mathbb{E}\left\|\int_{0}^{t}G_h(t-s)\phi(s) dB^H(s)\right\|^2\nonumber\\
&\leq& 2\mathbb{E}\left\|\int_{0}^{t}G_h(t-s)\left(\phi(t)-\phi(s)\right) dB^H(s)\right\|^2\nonumber\\
&+&2\mathbb{E}\left\|\int_{0}^{t}G_h(t-s)\phi(t) dB^H(s)\right\|^2\nonumber\\
&\leq& 2C\int_{0}^{t}\left\|G_h(t-s)\left(\phi(t)-\phi(s)\right)\right\|^2_{L^0_2} ds\nonumber\\
&+&2C_H\sum_{i\in\N^d}\left(\int_0^t\left\|G_h(t-s)\phi(t) Q^{\frac 12}e_i\right\|^{\frac 1H}ds\right)^{2H}\nonumber\\
&\leq& C h^{2(2H+\beta-1)}\int_{0}^{t}(t-s)^{-2H}\left\|A^{\frac{\beta-1}2}\left(\phi(t)-\phi(s)\right)\right\|^2_{L^0_2} ds\nonumber\\
&&+C_H\sum_{i\in\N^d}Ch^{2(2H+\beta-1)}\left\|A^{\frac{\beta-1}2}\phi(t) Q^{\frac 12}e_i\right\|^2\nonumber\\
&\leq& 
C h^{2(2H+\beta-1)}\int_{0}^{t}(t-s)^{2\delta-2H} ds+C h^{2(2H+\beta-1)}\left\|A^{\frac{\beta-1}2}\phi(t)\right\|^2_{L^0_2}\nonumber\\
&\leq& 
C h^{2(2H+\beta-1)}t^{2\delta-2H+1}+C h^{2(2H+\beta-1)}\nonumber\\
&\leq& C h^{2(2H+\beta-1)}.
\end{eqnarray}
Combining the estimates \eqref{space_err0}, \eqref{space_err11}, \eqref{spa_err_2} and applying  Gronwall inequality  ends the proof. $\hfill\square$

\subsection{Proof of Theorem \ref{cstrongconvthm} for implicit scheme}
It is important to mention that the estimates made in this section are inspired by the results in \cite[(4.7)-(4.14), (4.25)-(4.29)]{Wanb},
 when the linear operator$ A $ is self-adjoint.  For our case where A is not necessarily self-adjoint,  let us present  some preparatory results.
\begin{lem}
\label{lemim1}
For any $m$, $h$ and $\Delta t$ the following estimates holds
\begin{enumerate}
\item[(i)]
\begin{eqnarray}
\label{lemimeq1}
\left\|(I+\Delta t A_h)^{-m}\right\|\leq 1.
\end{eqnarray}
\item[(ii)] For all $u\in D(A^{\frac{\gamma-1}2})$, $0\leq \gamma\leq 1$,
\begin{eqnarray}
\label{lemimeq2}
\left\|(S^m_{h,\Delta t}-S_h(t_m))P_hu\right\|\leq C\Delta t^{\frac{2H+\gamma-1}2}t_m^{-H}\|u\|_{\gamma-1}.
\end{eqnarray}
\item[(iii)] For all $u\in D(A^{\frac{\gamma-1}2})$, $0\leq \gamma\leq 1$,
\begin{eqnarray}
\label{lemimeq3}
\left\|(S^m_{h,\Delta t}-S_h(t))P_hu\right\|\leq C\Delta t^{\frac{2H+\gamma-1}2}t^{-H}\|u\|_{\gamma-1},\hspace{0.5cm}t\in[t_{j-1},t_j),
\end{eqnarray}
for any $j=1,2,\cdot\cdot\cdot,M$.
\item[(iv)] If $u\in D(A^{\mu/2})$, $0\leq \mu\leq 2$, then
\begin{eqnarray}
\label{lemimeq4}
\|(S^m_{h,\Delta t}-S_h(t_m))P_hu\|\leq C\Delta t^{\mu/2}\|u\|_{\mu}.
\end{eqnarray}
\item [(v)] For all non smooth data $u\in\mathcal{H}$,
\begin{eqnarray}
\label{lemimeq5}
\|(S^m_{h,\Delta t}-S_h(t_m))P_hu\|\leq C\Delta t\hspace{0.1cm}t_m^{-1}\|u\|.
\end{eqnarray}
\end{enumerate}
\end{lem}
\textit{Proof} see \cite[Lemma 3.3]{Mukb} for the proof of $(i)$, $(iv)$ and $(v)$.   For the proof of $(ii)$,
 we use\cite[Lemma 3.3 $(iii)$ (88)]{Mukb} as follows
\begin{eqnarray}
\label{proof}
\left\|K_h(m)A_h^{\frac{1-\gamma}2}\right\|_{L(\mathcal{H})}&\leq& C t_m^{-3/2+\gamma/2}\Delta t=C t_m^{-H}t_m^{\frac{2H+\gamma-3}2}\Delta t\nonumber\\
&\leq& C t_m^{-H}t_1^{\frac{2H+\gamma-3}2}\Delta t=C t_m^{-H}\Delta t^{\frac{2H+\gamma-1}2}.
\end{eqnarray}
Hence substituting \eqref{proof} in \cite[(84)]{Mukb} completes the proof of $(ii)$. Now for the proof of $(iii)$, triangle inequality, Lemma \ref{lemim1} $(ii)$, the property of discrete semigroup and \cite[(83)]{Mukb} yields
\begin{eqnarray}
\label{proof2}
&&\left\|(S^m_{h,\Delta t}-S_h(t))P_hu\right\|\nonumber\\
&\leq&\left\|(S^m_{h,\Delta t}-S_h(t_j))P_hu\right\|+\left\|(S_h(t_j)-S_h(t))P_hu\right\|\nonumber\\
&\leq& \left\|S_h(t)A_h^H\right\|_{L(\mathcal{H})}\left\|A_h^{-\frac{2H+\gamma-1}2}(S_h(t_j-t)-I)\right\|_{L(\mathcal{H})}\left\|A_h^{\frac{\gamma-1}2}P_hu\right\|\nonumber\\
&&+C\Delta t^{\frac{2H+\gamma-1}2}t_j^{-H}\|u\|_{\gamma-1}\nonumber\\
&\leq& C t^{-H}(t_j-t)^{\frac{2H+\gamma-1}2}\left\|A^{\frac{\gamma-1}2}u\right\|+C\Delta t^{\frac{2H+\gamma-1}2}t_j^{-H}\|u\|_{\gamma-1}\nonumber\\
&\leq& C\Delta t^{\frac{2H+\gamma-1}2}t^{-H}\|u\|_{\gamma-1}.
\end{eqnarray}

The proof of Lemma \ref{lemim1} is thus completed. $\hfill{\square}$
\begin{lem} 
\label{lemim2}
\begin{enumerate}
\item[(i)] For any $\rho\in[0,1]$ and $u\in D(A^{-\rho/2})$ there exists a positive constant $C$ such that 
\begin{eqnarray}
\label{lemimeq6}
\left\|\sum_{j=1}^{m}\int_{t_{j-1}}^{t_j}(S^j_{h,\Delta t}-S_h(s))P_huds\right\|\leq C\Delta t^{\frac{2-\rho}2-\epsilon}\|u\|_{-\rho}.
\end{eqnarray}
\item[(ii)] For any $\mu\in[0,1]$ and $u\in D(A^{\frac{\mu-1}2})$ the following estimate holds
\begin{eqnarray}
\label{lemimeq7}
\left(\sum_{j=1}^{m}\int_{t_{j-1}}^{t_j}\|(S^j_{h,\Delta t}-S_h(s))P_hu\|^{\frac 1H}ds\right)^{2H}\leq C\Delta t^{2H+\mu-1-\epsilon}\|u\|_{\mu-1}^2.
\end{eqnarray}
\end{enumerate}
Where $\epsilon$ is an arbitrary small number.
\end{lem}
\textit{Proof} See \cite[Lemma 3.5]{Mukb} for the proof of $(i)$. For the proof of $(ii)$, we have
\begin{eqnarray}
\label{lemim_eq7}
&&\sum_{j=1}^{m}\int_{t_{j-1}}^{t_j}\|(S^j_{h,\Delta t}-S_h(s))P_hu\|^{\frac 1H}ds\nonumber\\
&=&\int_{0}^{\Delta t}\|(S^1_{h,\Delta t}-S_h(s))P_hu\|^{\frac 1H}ds+\sum_{j=2}^{m}\int_{t_{j-1}}^{t_j}\|(S^j_{h,\Delta t}-S_h(s))P_hu\|^{\frac 1H}ds\nonumber\\
&=:& K_1+K_2.
\end{eqnarray} 
Using triangle inequality and the estimate $(a+b)^k\leq 2^{k-1}a^k+2^{k-1}b^k$ (with $k\geq1$ and $a,b\geq0$) we obtain
\begin{eqnarray}
\label{lemim_eq71}
K_1&\leq& 2^{\frac 1H-1}\int_{0}^{\Delta t}\|(S^1_{h,\Delta t}-S_h(t_1))P_hu\|^{\frac 1H}ds\nonumber\\
&+&2^{\frac 1H-1}\int_{0}^{\Delta t}\|(S_h(t_1)-S_h(s))P_hu\|^{\frac 1H}ds\nonumber\\
&=:&2^{\frac 1H-1}K_{11}+2^{\frac 1H-1}K_{12}.
\end{eqnarray}
By Lemma \ref{lemim1} $(ii)$ with $\gamma=\mu$ we obtain
\begin{eqnarray}
\label{lemim_eq711}
K_{11}&\leq& C\int_{0}^{\Delta t}t_1^{-1}\Delta t^{\frac{2H+\mu-1}{2H}}\|u\|^{\frac 1H}_{\mu-1}ds\nonumber\\
&\leq& C\Delta t^{\frac{2H+\mu-1}{2H}}\|u\|^{\frac 1H}_{\mu-1}.
\end{eqnarray}
By inserting an appropriate power of $A_h$, \cite[(81)]{Mukb} and Remark \ref{disc_semigroup} yields
\begin{eqnarray}
\label{lemim_eq712}
&&K_{12}\nonumber\\
&\leq& \int_{0}^{\Delta t}\left\|\left(S_h(t_1-s)-I\right)S_h(s)A_h^{\frac{1-\mu}2}\right\|_{L(\mathcal{H})}^{\frac 1H}\left\|A_h^{\frac{\mu-1}2}P_hu\right\|^{\frac 1H}ds\nonumber\\
&\leq& C\int_{0}^{\Delta t}\left\|\left(S_h(t_1-s)-I\right)A_h^{\frac{-2H-\mu+1}2}\right\|_{L(\mathcal{H})}^{\frac 1H}\left\|S_h(s)A_h^H\right\|_{L(\mathcal{H})}^{\frac 1H}\left\|A_h^{\frac{\mu-1}2}P_hu\right\|^{\frac 1H}ds\nonumber\\
&\leq& C\int_{0}^{\Delta t}(t_1-s)^{\frac{2H+\mu-1}{2H}}\left\|S_h(s)A_h^H\right\|_{L(\mathcal{H})}^{\frac 1H}\left\|A_h^{\frac{\mu-1}2}u\right\|^{\frac 1H}ds\nonumber\\
&\leq& C\Delta t^{\frac{2H+\mu-1}{2H}}\left\|u\right\|^{\frac 1H}_{\mu-1}\left(\int_{0}^{\Delta t}\left\|S_h(s)A_h^H\right\|_{L(\mathcal{H})}^{\frac 1H}ds\right)\nonumber\\
&\leq& C\Delta t^{\frac{2H+\mu-1}{2H}}\left\|u\right\|^{\frac 1H}_{\mu-1}\left(\int_{0}^{\Delta t}\left\|S_h(\Delta t-s)A_h^H\right\|_{L(\mathcal{H})}^{\frac 1H}ds\right)\nonumber\\
&\leq& C\Delta t^{\frac{2H+\mu-1}{2H}}\left\|u\right\|^{\frac 1H}_{\mu-1}.
\end{eqnarray} 
Substituting \eqref{lemim_eq711} and \eqref{lemim_eq712} in \eqref{lemim_eq71} yields
\begin{eqnarray}
\label{lemim_eq7_1}
K_1\leq C\Delta t^{\frac{2H+\mu-1}{2H}}\left\|u\right\|^{\frac 1H}_{\mu-1}.
\end{eqnarray}
Concerning the estimate of $K_2$, let $\epsilon>0$ small enough, Lemma \ref{lemim1} $(iii)$ with $\gamma=\mu$ yields
\begin{eqnarray}
\label{lemim_eq72}
K_2&=&\sum_{j=2}^{m}\int_{t_{j-1}}^{t_j}\|(S^j_{h,\Delta t}-S_h(s))P_hu\|^{\frac 1H}ds\nonumber\\
&\leq& C\sum_{j=2}^{m}\int_{t_{j-1}}^{t_j}\Delta t^{\frac{2H+\mu-1}{2H}}s^{-1}\|u\|^{\frac 1H}_{\mu-1}ds\nonumber\\
&\leq& C\left(\sum_{j=2}^{m}\int_{t_{j-1}}^{t_j}s^{-1+\frac{\epsilon}{2H}}s^{-\frac{\epsilon}{2H}}ds\right)\Delta t^{\frac{2H+\mu-1}{2H}}\|u\|^{\frac 1H}_{\mu-1}\nonumber\\
&\leq& C\left(\sum_{j=2}^{m}\int_{t_{j-1}}^{t_j}s^{-1+\frac{\epsilon}{2H}}t_1^{-\frac{\epsilon}{2H}}ds\right)\Delta t^{\frac{2H+\mu-1}{2H}}\|u\|^{\frac 1H}_{\mu-1}\nonumber\\
&\leq& C\left(\sum_{j=2}^{m}\int_{t_{j-1}}^{t_j}s^{-1+\frac{\epsilon}{2H}}ds\right)\Delta t^{\frac{2H+\mu-1-\epsilon}{2H}}\|u\|^{\frac 1H}_{\mu-1}\nonumber\\
&\leq& C\Delta t^{\frac{2H+\mu-1-\epsilon}{2H}}\|u\|^{\frac 1H}_{\mu-1}.
\end{eqnarray} 
Adding \eqref{lemim_eq7_1} and \eqref{lemim_eq72} yields
\begin{eqnarray}
\label{lemimeq72}
\sum_{j=1}^{m}\int_{t_{j-1}}^{t_j}\|(S^j_{h,\Delta t}-S_h(s))P_hu\|^{\frac 1H}ds\leq C\Delta t^{\frac{2H+\mu-1-\epsilon}{2H}}\|u\|^{\frac 1H}_{\mu-1},
\end{eqnarray}
hence
\begin{eqnarray}
\label{lemimeq7*}
\left(\sum_{j=1}^{m}\int_{t_{j-1}}^{t_j}\|(S^j_{h,\Delta t}-S_h(s))P_hu\|^{\frac 1H}ds\right)^{2H}\leq C\Delta t^{2H+\mu-1-\epsilon}\|u\|^2_{\mu-1}.
\end{eqnarray}
This completes the proof of Lemma \ref{lemim2}. $\hfill\square$
 
With these two lemmas, we are  now ready to  prove our theorem for the implicit scheme . 
 In fact, using the standard technique in the error analysis, we split the fully discrete error in two terms as
\begin{eqnarray*}
\|X(t_m)-X^h_m\|_{L^2(\Omega,\mathcal{H})}&\leq& \|X(t_m)-X^h(t_m)\|_{L^2(\Omega,\mathcal{H})}+\|X^h(t_m)-X^h_m\|_{L^2(\Omega,\mathcal{H})}\\
&=:&err_0+err_1.
\end{eqnarray*}
Note that the space error $err_0$ is estimated by Lemma \ref{serror}. It remains to estimate the time error $err_1$.\\
We recall that the exact solution at $t_m$ of the semidiscrete problem \eqref{pb2} is given by
\begin{eqnarray}
\label{mild}
X^h(t_m)&=&S_h(t_m)X^h_0+\int_{0}^{t_m}S_h(t_m-s)P_hF(X^h(s))ds\nonumber\\
&+&\int_{0}^{t_m}S_h(t_m-s)P_h\phi(s) dB^H(s).
\end{eqnarray}
We also recall that the numerical solution at $t_m$ given by \eqref{impl} can be rewritten as 
\begin{eqnarray}
\label{impl1}
X^h_m&=&S_{h,\Delta t}^mX^h_0 +\int_{0}^{t_m}S_{h,\Delta t}^{m-[s]^m}P_hF(X^h_{[s]^m})ds\nonumber\\
&+&\int_{0}^{t_m}S_{h,\Delta t}^{m-[s]^m}P_h\phi([s]) dB^H(s),
\end{eqnarray}
where the notation $[t]$, $[t]^m$  are defined by
\begin{eqnarray}
\label{notation}
[t]:=\left[\frac{t}{\Delta t}\right]\Delta t\hspace{1cm}\text{and}\hspace{1cm}[t]^m:=\left[\frac{t}{\Delta t}\right]
\end{eqnarray}
It follows from \eqref{mild} and \eqref{impl1} that
\begin{eqnarray}
\label{err1}
err_1&\leq& \left\|(S_h(t_m)-S^m_{h,\Delta t})P_hX_0\right\|_{L^2(\Omega,\mathcal{H})}\nonumber\\
&&+\left\|\int_{0}^{t_m}S_h(t_m-s)P_hF(X^h(s))-S_{h,\Delta t}^{m-[s]^m}P_hF(X^h_{[s]^m})ds\right\|_{L^2(\Omega,\mathcal{H})}\nonumber\\
&&+\left\|\int_{0}^{t_m}S_h(t_m-s)P_h\phi(s)-S_{h,\Delta t}^{m-[s]^m}P_h\phi([s]) dB^H(s)\right\|_{L^2(\Omega,\mathcal{H})}\nonumber\\
&=:&I_0+I_1+I_2.
\end{eqnarray}
As we said at the beginning of this section, following closely the work done in \cite[(4.7)-(4.14)]{Wanb} and  replacing its preparatory results with Lemma \ref{lemim1} $(i)$, $(iv)$ with $\mu=2H+\beta-1$, $(v)$, Lemma \ref{lemim2} $(i)$ with $\rho=0$, 
Remark \ref{disc_semigroup} \eqref{semigroup_prp1} with $\delta=\gamma=1$, Assumptions \ref{non}-\ref{init}, boundedness of $S_h(t_m-s)$ and $P_h$, the stability properties of a discrete semigroup $S_h(t)$, \eqref{ass_non_lin1} and \eqref{reg2}, we have
\begin{eqnarray}
\label{det_err1}
I_0+I_1\leq C\Delta t^{\frac{2H+\beta-1}2}+C\Delta t^{1-\epsilon}+C\Delta t\sum_{i=0}^{m-1}\|X^h(t_i)-X^h_i\|_{L^2(\Omega,\mathcal{H})}.
\end{eqnarray}
Note that in this work, we do not need to impose an assumption  on $F''$ to increase the convergence rate as it is done in \cite{Wanb}. Indeed, thanks to \eqref{reg2} the following  estimate is largely sufficient to reach a higher rate.

\begin{eqnarray}
\label{deterr}
I_{11}&\leq&\sum_{i=0}^{m-1} \int_{t_i}^{t_{i+1}}\|S_h(t_m-s)P_h(F(X^h(s))-F(X^h(t_i)))\|_{L^2(\Omega,\mathcal{H})}ds\nonumber\\
&\leq& C\sum_{i=0}^{m-1} \int_{t_i}^{t_{i+1}}\|X^h(s)-X^h(t_i)\|_{L^2(\Omega,\mathcal{H})}ds\nonumber\\
&\leq& C\Delta t^{\frac{2H+\beta-1}2}.
\end{eqnarray}

Let us focus  now on the estimate $I_2$, using triangle inequality and the estimate $(a+b)^2\leq 2a^2+2b^2$ we split it in three terms
\begin{eqnarray}
\label{stoerr}
I_2^2&=&\left\|\int_{0}^{t_m}S_h(t_m-s)P_h\phi(s)-S_{h,\Delta t}^{m-[s]^m}P_h\phi([s]) dB^H(s)\right\|^2_{L^2(\Omega,\mathcal{H})}\nonumber\\
&\leq& 2\left\|\int_{0}^{t_m}S_h(t_m-s)P_h\left(\phi(s)-\phi([s])\right)dB^H(s)\right\|^2_{L^2(\Omega,\mathcal{H})}\nonumber\\
&&+2\left\|\int_{0}^{t_m}\left(S_h(t_m-s)-S_{h,\Delta t}^{m-[s]^m}\right)P_h\phi([s]) dB^H(s)\right\|^2_{L^2(\Omega,\mathcal{H})}\nonumber\\
&\leq& 2\left\|\int_{0}^{t_m}S_h(t_m-s)P_h\left(\phi(s)-\phi([s])\right)dB^H(s)\right\|^2_{L^2(\Omega,\mathcal{H})}\nonumber\\
&&+2\left\|\int_{0}^{t_m}\left(S_h(t_m-s)-S_{h,\Delta t}^{m-[s]^m}\right)P_h\left(\phi([s])-\phi(t_{m-1})\right) dB^H(s)\right\|^2_{L^2(\Omega,\mathcal{H})}\nonumber\\
&&+2\left\|\int_{0}^{t_m}\left(S_h(t_m-s)-S_{h,\Delta t}^{m-[s]^m}\right)P_h\phi(t_{m-1}) dB^H(s)\right\|^2_{L^2(\Omega,\mathcal{H})}\nonumber\\&=:& 2I_{21}^2+4I_{22}^2+4I_{23}^2.
\end{eqnarray}
Firstly using \eqref{bound}, inserting an appropriate power of $A_h$, \cite[(81)]{Mukb}, \assref{noise} (more precisely \eqref{ass_noise_term3}) and Remark \ref{disc_semigroup} ( \eqref{semigroup_prp3} with $\rho=1-\beta$ ) we obtain
\begin{eqnarray}
\label{stoerr1}
I_{21}^2&=&\left\|\int_{0}^{t_m}S_h(t_m-s)P_h\left(\phi(s)-\phi([s])\right)dB^H(s)\right\|^2_{L^2(\Omega,\mathcal{H})}\nonumber\\
&\leq& C\int_{0}^{t_m}\left\|S_h(t_m-s)P_h\left(\phi(s)-\phi([s])\right)\right\|^2_{L^0_2}ds\nonumber\\
&\leq& C\int_{0}^{t_m}\left\|S_h(t_m-s)A_h^{\frac{1-\beta}2}\right\|^2_{L(\mathcal{H})}\left\|A_h^{\frac{\beta-1}2}P_h\left(\phi(s)-\phi([s])\right)\right\|^2_{L^0_2}ds\nonumber\\
&\leq& C\int_{0}^{t_m}\left\|S_h(t_m-s)A_h^{\frac{1-\beta}2}\right\|^2_{L(\mathcal{H})}\left\|A^{\frac{\beta-1}2}\left(\phi(s)-\phi([s])\right)\right\|^2_{L^0_2}ds\nonumber\\
&\leq& C\int_{0}^{t_m}(s-[s])^{2\delta}\left\|S_h(t_m-s)A_h^{\frac{1-\beta}2}\right\|^2_{L(\mathcal{H})}ds\nonumber\\
&\leq& C\Delta t^{2\delta}\int_{0}^{t_m}\left\|S_h(t_m-s)A_h^{\frac{1-\beta}2}\right\|^2_{L(\mathcal{H})}ds\nonumber\\
&\leq& C\Delta t^{2\delta}t_m^{\beta}\leq C\Delta t^{2H+\beta-1}.
\end{eqnarray}
Secondly \eqref{bound}, the change of variable $j=m-k$ and $\varsigma=t_m-s$, \lemref{lemim1} $(iii)$  with $\gamma=\beta$ and \assref{noise} ( more precisely \eqref{ass_noise_term3}) yields 
\begin{eqnarray}
\label{stoerr2}
I_{22}^2&=&\left\|\int_{0}^{t_m}\left(S_h(t_m-s)-S_{h,\Delta t}^{m-[s]^m}\right)P_h\left(\phi([s])-\phi(t_{m-1})\right) dB^H(s)\right\|^2_{L^2(\Omega,\mathcal{H})}\nonumber\\
&\leq&C\int_{0}^{t_m}\left\|\left(S_h(t_m-s)-S_{h,\Delta t}^{m-[s]^m}\right)P_h\left(\phi([s])-\phi(t_{m-1})\right)\right\|^2_{L^0_2} ds\nonumber\\
&\leq&C\sum_{k=0}^{m-1}\int_{t_k}^{t_{k+1}}\left\|\left(S_h(t_m-s)-S_{h,\Delta t}^{m-k}\right)P_h\left(\phi(t_k)-\phi(t_{m-1})\right)\right\|^2_{L^0_2} ds\nonumber\\
&\leq&C\sum_{j=1}^{m}\int_{t_{j-1}}^{t_j}\left\|\left(S_h(\varsigma)-S_{h,\Delta t}^{j}\right)P_h\left(\phi(t_{m-j})-\phi(t_{m-1})\right)\right\|^2_{L^0_2} d\varsigma\nonumber\\
&\leq&C\sum_{j=1}^{m}\int_{t_{j-1}}^{t_j}\Delta t^{2H+\beta-1}\varsigma^{-2H}\left\|A^{\frac{\beta-1}2}\left(\phi(t_{m-j})-\phi(t_{m-1})\right)\right\|^2_{L^0_2} d\varsigma\nonumber\\
&\leq&C\Delta t^{2H+\beta-1}\left(\sum_{j=1}^{m}\int_{t_{j-1}}^{t_j}\varsigma^{-2H}t_{j-1}^{2\delta} d\varsigma\right)\nonumber\\
&\leq&C\Delta t^{2H+\beta-1}\left(\sum_{j=1}^{m}\int_{t_{j-1}}^{t_j}\varsigma^{2\delta-2H} d\varsigma\right)\nonumber\\
&\leq&C\Delta t^{2H+\beta-1}\left(\int_{0}^{t_m}\varsigma^{2\delta-2H} d\varsigma\right)\nonumber\\
&\leq&C\Delta t^{2H+\beta-1}t_m^{2\delta-2H+1}\leq C\Delta t^{2H+\beta-1}.
\end{eqnarray}
Thirdly using \eqref{stoint1}, the change of variable $j=m-k$ and $\varsigma=t_m-s$, \lemref{lemim2} $(ii)$  with $\mu=\beta$ and \assref{noise} ( more precisely \eqref{ass_noise_term2} ) we obtain
\begin{eqnarray}
\label{stoerr3}
I_{23}^2&=& \left\|\int_{0}^{t_m}\left(S_h(t_m-s)-S_{h,\Delta t}^{m-[s]^m}\right)P_h\phi(t_{m-1}) dB^{H}(s)\right\|^2_{L^2(\Omega,\mathcal{H})}\nonumber\\
&\leq& C_H\sum_{i\in\N^d}\left(\int_{0}^{t_m}\left\|\left(S_h(t_m-s)-S_{h,\Delta t}^{m-[s]^m}\right)P_h\phi(t_{m-1}) Q^{\frac 12}e_i\right\|^{\frac 1H}ds\right)^{2H}\nonumber\\
&\leq& C_H\sum_{i\in\N^d}\left(\sum_{k=0}^{m-1}\int_{t_k}^{t_{k+1}}\left\|\left(S_h(t_m-s)-S_{h,\Delta t}^{m-k}\right)P_h\phi(t_{m-1}) Q^{\frac 12}e_i\right\|^{\frac 1H}ds\right)^{2H}\nonumber\\
&\leq& C_H\sum_{i\in\N^d}\left(\sum_{j=1}^{m}\int_{t_{j-1}}^{t_j}\left\|\left(S_h(\varsigma)-S_{h,\Delta t}^j\right)P_h\phi(t_{m-1}) Q^{\frac 12}e_i\right\|^{\frac 1H}d\varsigma\right)^{2H}\nonumber\\
&\leq& C_H\sum_{i\in\N^d}\left(\Delta t^{2H+\beta-1-\epsilon}\|\phi(t_{m-1}) Q^{\frac 12}e_i\|^2_{\beta-1}\right)\nonumber\\
&\leq& C\Delta t^{2H+\beta-1-\epsilon}\|A^{\frac{\beta-1}2}\phi(t_{m-1})\|^2_{L^0_2}\leq C\Delta t^{2H+\beta-1-\epsilon}.
\end{eqnarray}
Hence inserting \eqref{stoerr1}-\eqref{stoerr3} in \eqref{stoerr} and taking the square-root gives 
\begin{eqnarray}
\label{sto_err1}
I_2\leq C\Delta t^{\frac{2H+\beta-1}2-\epsilon}.
\end{eqnarray} 
Adding \eqref{det_err1} and \eqref{sto_err1} we obtain
\begin{eqnarray*}
err_1&=& \|X^h(t_m)-X^h_m\|_{L^2(\Omega,\mathcal{H})}\nonumber\\
&\leq& C\Delta t^{\frac{2H+\beta-1}2}+C\Delta t^{1-\varepsilon}+C\Delta t^{\frac{2H+\beta-1}2-\varepsilon}\nonumber\\
&+&C\Delta t\sum_{i=0}^{m-1}\|X^h(t_i)-X^h_i\|_{L^2(\Omega,\mathcal{H})}.
\end{eqnarray*}
Applying the discrete version of the Gronwall inequality yields
\begin{eqnarray}
\label{time_err1}
err_1\leq C\Delta t^{\frac{2H+\beta-1}2-\epsilon}.
\end{eqnarray}
Adding \eqref{spa_err} and \eqref{time_err1} completes the proof.$\hfill\square$

In what follows, we will present a corollary of Theorem \ref{cstrongconvthm} for the implicit Euler scheme where the linear operator $A$ is assumed to be self-adjoint.  The optimal strong convergence rate in time  $\mathcal{O}(\Delta t)$ is reached.
\begin{cor}
\label{impl1self-adj}
Let $X(t_m)$ be the mild solution of \eqref{pb1} ($A$ self-adjoint) at time $t_m=m\Delta t$, $\Delta t\geq 0$ represented by \eqref{mild_sol_pb1}. Let $X^h_m$ be the numerical approximation through \eqref{impl}. Under Assumptions \ref{noise}-\ref{init}, $\beta\in(0,1]$, then the following holds
\begin{eqnarray}
\label{strongconvimpl}
\left(\mathbb{E}\Vert X(t_m)-X^h_m\Vert^2\right)^{\frac 12}\leq C\left(h^{2H+\beta-1}+\Delta t^{\frac{2H+\beta-1}2}\right)
\end{eqnarray}
\end{cor}
For the proof of this corollary, we need to update our preparatory results, more precisely \lemref{lemim2} in the self-adjoint case. 
The result is presented in the following lemma:
\begin{lem}
\label{lemim2*}
\begin{enumerate}
\item[(i)] For any $\rho\in[0,1]$ and $u\in D(A^{-\rho/2})$ there exists a positive constant $C$ such that 
\begin{eqnarray}
\label{lemimeq6*}
\left\|\sum_{j=1}^{m}\int_{t_{j-1}}^{t_j}(S^j_{h,\Delta t}-S_h(s))P_huds\right\|\leq C\Delta t^{\frac{2-\rho}2}\|u\|_{-\rho}.
\end{eqnarray}
\item[(ii)] For any $\mu\in[0,1]$ and $x\in D(A^{\frac{\mu-1}2})$ the following estimate holds
\begin{eqnarray}
\label{lemimeq7'}
&&\sum_{i,j=0}^{m-1}\int_{t_j}^{t_{j+1}}\int_{t_i}^{t_{i+1}}\left\langle\mathcal{S}_h(u,t_i)P_hx,\mathcal{S}_h(v,t_j)P_hx\right\rangle \kappa(u,v)dudv\nonumber\\
&\leq& C\Delta t^{2H+\mu-1}\|x\|_{\mu-1}^2,
\end{eqnarray}
and
\begin{eqnarray}
\label{lemimeq7"}
&&\sum_{i,j=0}^{m-1}\int_{t_j}^{t_{j+1}}\int_{t_i}^{t_{i+1}}\left\langle\mathcal{T}_h(i)P_hx,\mathcal{T}_h(j)P_hx\right\rangle \kappa(u,v)dudv\nonumber\\
&\leq& C\Delta t^{2H+\mu-1}\|x\|_{\mu-1}^2,
\end{eqnarray}
\end{enumerate}
where
\begin{eqnarray}
\label{not}
\begin{split}
\mathcal{S}_h(u,t_i):=S_h(t_m-u)-S_h(t_m-t_i),\quad \mathcal{T}_h(i):=S_h(t_m-t_i)-S^{m-i}_{h,\Delta t},\\
\text{and}\quad \kappa(u,v):=H(2H-1)\vert u-v\vert^{2H-2}.
\end{split}
\end{eqnarray}
\end{lem}
\textit{Proof.} See \cite[Proof of Lemma 4.4 $(i)$]{Kru} for the proof of $(i)$ and \cite[Lemmas 4.8 and 4.9]{Wanc}, \cite[(83)]{Mukb} for the proof of $(ii)$. $\hfill\square$\\

With this new lemma, we are now in position to prove our \coref{impl1self-adj}. 

\textit{Proof of \coref{impl1self-adj}.}
Recall that the time error $err_1$ is defined as
\begin{eqnarray}
\label{err1*}
err_1&\leq& \left\|(S_h(t_m)-S^m_{h,\Delta t})P_hX_0\right\|_{L^2(\Omega,\mathcal{H})}\nonumber\\
&+&\left\|\int_{0}^{t_m}S_h(t_m-s)P_hF(X^h(s))-S_{h,\Delta t}^{m-[s]^m}P_hF(X^h_{[s]^m})ds\right\|_{L^2(\Omega,\mathcal{H})}\nonumber\\
&+&\left\|\int_{0}^{t_m}S_h(t_m-s)P_h\phi(s)-S_{h,\Delta t}^{m-[s]^m}P_h\phi([s]) dB^H(s)\right\|_{L^2(\Omega,\mathcal{H})}\nonumber\\
&=:&I_0+I_1+I_2.
\end{eqnarray}
Following closely the work done in \cite[(4.7)-(4.14)]{Wanb} and  replacing its preparatory results with Lemma \ref{lemim1} $(i)$, $(iv)$ with $\mu=2H+\beta-1$, $(v)$ with $\sigma=0$, Lemma \ref{lemim2*} $(i)$ with $\rho=0$, Remark \ref{disc_semigroup} \eqref{semigroup_prp1} with $\delta=\gamma=1$, Assumptions \ref{non}-\ref{init}, boundedness of $S_h(t_m-s)$ and $P_h$, the stability properties of a discrete semigroup $S_h(t)$, \eqref{ass_non_lin1} and \eqref{reg2}, we have
\begin{eqnarray}
\label{det_err1*}
I_0+I_1\leq C\Delta t^{\frac{2H+\beta-1}2}+C\Delta t+C\Delta t\sum_{i=0}^{m-1}\|X^h(t_i)-X^h_i\|_{L^2(\Omega,\mathcal{H})}.
\end{eqnarray} 
Concerning the estimate $I_2$, we also split it in three terms as in \eqref{stoerr}. The estimates $I_{21}^2$ and $I_{22}^2$ still the same  but  
we need to re-estimate $I_{23}^2$. In this fact, since the sequence of random variables 

$\left(\int_0^{t_m}\left(S_h(t_m-s)-S_{h,\Delta t}^{m-[s]^m}\right)
P_h\phi(t_{m-1}) Q^{\frac 12}e_id\beta_i^H(s),\,\,i\in\N^d\right)$ are mutually independent Gaussian random variable, using the estimate $(a+b)^2\leq 2a^2+2b^2$, notation \eqref{not}, Assumption \ref{noise} ( more precisely \eqref{ass_noise_term2}),  \lemref{lemim2*} $(ii)$  with $\mu=\beta$, we obtain

{\small
\begin{eqnarray}
\label{stoerr3*}
I_{23}^2&=& \left\|\int_{0}^{t_m}\left(S_h(t_m-s)-S_{h,\Delta t}^{m-[s]^m}\right)P_h\phi(t_{m-1}) dB^{H}(s)\right\|^2_{L^2(\Omega,\mathcal{H})}\nonumber\\
&=&\left\|\sum_{i\in\N^d}\int_{0}^{t_m}\left(S_h(t_m-s)-S_{h,\Delta t}^{m-[s]^m}\right)P_h\phi(t_{m-1}) Q^{\frac 12}e_id\beta_i^H(s)\right\|^2_{L^2(\Omega,\mathcal{H})}\nonumber\\
&=&\sum_{i\in\N^d}\left\|\int_{0}^{t_m}\left(S_h(t_m-s)-S_{h,\Delta t}^{m-[s]^m}\right)P_h\phi(t_{m-1}) Q^{\frac 12}e_id\beta_i^H(s)\right\|^2_{L^2(\Omega,\mathcal{H})}\nonumber\\
&\leq& 2\sum_{i\in\N^d}\left\|\int_{0}^{t_m}\left(S_h(t_m-s)-S_h(t_m-[s])\right)P_h\phi(t_{m-1}) Q^{\frac 12}e_id\beta_i^H(s)\right\|^2_{L^2(\Omega,\mathcal{H})}\nonumber\\
&&+2\sum_{i\in\N^d}\left\|\int_{0}^{t_m}\left(S_h(t_m-[s])-S_{h,\Delta t}^{m-[s]^m}\right)P_h\phi(t_{m-1}) Q^{\frac 12}e_id\beta_i^H(s)\right\|^2_{L^2(\Omega,\mathcal{H})}\nonumber\\
&=& 2\sum_{i\in\N^d}\left\|\int_{0}^{t_m}\mathcal{S}_h(s,[s])P_h\phi(t_{m-1}) Q^{\frac 12}e_id\beta_i^H(s)\right\|^2_{L^2(\Omega,\mathcal{H})}\nonumber\\
&&+2\sum_{i\in\N^d}\left\|\int_{0}^{t_m}\mathcal{T}_h([s]^m)P_h\phi(t_{m-1}) Q^{\frac 12}e_id\beta_i^H(s)\right\|^2_{L^2(\Omega,\mathcal{H})}\nonumber\\
&=&2\sum_{i\in\N^d}\int_{0}^{t_m}\int_{0}^{t_m}\left\langle \mathcal{S}_h(u,[u])P_h\phi(t_{m-1}) Q^{\frac 12}e_i,\mathcal{S}_h(v,[v])P_h\phi(t_{m-1}) Q^{\frac 12}e_i\right\rangle \kappa(u,v)dudv\nonumber\\
&&+2\sum_{i\in\N^d}\int_{0}^{t_m}\int_{0}^{t_m}\left\langle \mathcal{T}_h([u]^m)P_h\phi(t_{m-1}) Q^{\frac 12}e_i,\mathcal{T}_h([v]^m)P_h\phi(t_{m-1}) Q^{\frac 12}e_i\right\rangle \kappa(u,v)dudv\nonumber\\
&=&2\sum_{i\in\N^d}\left(\sum_{i,j=0}^{m-1}\int_{t_j}^{t_{j+1}}\int_{t_i}^{t_{i+1}}\left\langle\mathcal{S}_h(u,t_i)P_h\phi(t_{m-1}) Q^{\frac 12}e_i,\mathcal{S}_h(v,t_j)P_h\phi(t_{m-1}) Q^{\frac 12}e_i\right\rangle \kappa(u,v)dudv\right)\nonumber\\
&&+2\sum_{i\in\N^d}\left(\sum_{i,j=0}^{m-1}\int_{t_j}^{t_{j+1}}\int_{t_i}^{t_{i+1}}\left\langle\mathcal{T}_h(i)P_h\phi(t_{m-1}) Q^{\frac 12}e_i,\mathcal{T}_h(j)P_h\phi(t_{m-1}) Q^{\frac 12}e_i\right\rangle \kappa(u,v)dudv\right)\nonumber\\
&\leq& C\sum_{i\in\N^d}\left(\Delta t^{2H+\beta-1}\|\phi(t_{m-1}) Q^{\frac 12}e_i\|^2_{\beta-1}\right)\nonumber\\
&\leq& C\Delta t^{2H+\beta-1}\|A^{\frac{\beta-1}2}\phi(t_{m-1})\|^2_{L^0_2}\leq C\Delta t^{2H+\beta-1},
\end{eqnarray}
}
hence inserting \eqref{stoerr1}, \eqref{stoerr2} and \eqref{stoerr3*} in \eqref{stoerr} and taking the square-root gives 
\begin{eqnarray}
\label{sto_err1*}
I_2\leq C\Delta t^{\frac{2H+\beta-1}2}.
\end{eqnarray} 
Adding \eqref{det_err1*} and \eqref{sto_err1*} we obtain
\begin{eqnarray*}
err_1&=& \|X^h(t_m)-X^h_m\|_{L^2(\Omega,\mathcal{H})}\nonumber\\
&\leq& C\Delta t^{\frac{2H+\beta-1}2}+C\Delta t+C\Delta t\sum_{i=0}^{m-1}\|X^h(t_i)-X^h_i\|_{L^2(\Omega,\mathcal{H})}.
\end{eqnarray*}
Applying the discrete version of the Gronwall inequality yields
\begin{equation}
\label{time_err1*}
err_1\leq C\Delta t^{\frac{2H+\beta-1}2}.
\end{equation}
Adding \eqref{spa_err} and \eqref{time_err1*} completes the proof.$\hfill\square$
\subsection{Proof of \thmref{cstrongconvthm} for SETD1}
 As usual, spliting the fully discrete error in two terms yields
\begin{eqnarray*}
\|X(t_m)-Y^h_m\|_{L^2(\Omega,\mathcal{H})}&\leq& \|X(t_m)-X^h(t_m)\|_{L^2(\Omega,\mathcal{H})}+\|X^h(t_m)-Y^h_m\|_{L^2(\Omega,\mathcal{H})}\nonumber\\
&=:& err_0+err_2.
\end{eqnarray*}
Since the space error $err_0$ has been estimated by Lemma \ref{serror}, we only need to estimate the time error $err_2$.
Remember that the exact solution at $t_m$ is given by
\begin{eqnarray}
\label{mild1}
X^h(t_m)&=&S_h(t_m)X^h_0+\int_{0}^{t_m}S_h(t_m-s)P_hF(X^h(s))ds\nonumber\\
&+&\int_{0}^{t_m}S_h(t_m-s)P_h\phi(s) dB^H(s)
\end{eqnarray}
and we recall that the numerical solution at $t_m$ given by \eqref{SETD1} can be rewritten as 
\begin{eqnarray}
\label{SETD11}
Y^h_m&=&S_h(t_m)X^h_0 +\int_{0}^{t_m}S_h(t_m-s)P_hF(X^h_{[s]^m})ds\nonumber\\
&+&\int_{0}^{t_m}S_h(t_m-[s])P_h\phi([s]) dB^H(s),
\end{eqnarray}
where the notations $[t]$ and $[t]^m$ are given by \eqref{notation}. By \eqref{mild1} and \eqref{SETD11}, we have
\begin{eqnarray}
\label{err2}
err_2&\leq& \left\|\int_{0}^{t_m}S_h(t_m-s)P_h(F(X^h(s))-F(Y^h_{[s]^m}))ds\right\|_{L^2(\Omega,\mathcal{H})}\nonumber\\
&+&\left\|\int_{0}^{t_m}S_h(t_m-s)P_h\phi(s)-S_h(t_m-[s])P_h\phi([s])dB^H(s)\right\|_{L^2(\Omega,\mathcal{H})}\nonumber\\
&=:&I_1'+I_2'.
\end{eqnarray}
Applying the triangle inequality yields
\begin{eqnarray}
\label{deterr2}
I_1'&\leq& \left\|\int_{0}^{t_m}S_h(t_m-s)P_h(F(X^h(s))-F(X^h([s])))ds\right\|_{L^2(\Omega,\mathcal{H})}\nonumber\\
&+&\left\|\int_{0}^{t_m}S_h(t_m-s)P_h(F(X^h([s]))-F(Y^h_{[s]^m}))ds\right\|_{L^2(\Omega,\mathcal{H})}\nonumber\\
&=:&I_{11}'+I_{12}'.
\end{eqnarray}
 Using the boundedness of $P_h$ and $S_h(t_m-s)$, Lemma \ref{lem2} and \eqref{reg2}, we easily have
\begin{eqnarray}
\label{det_err21}
I_{11}'\leq C\Delta t^{\frac{2H+\beta-1}2},
\end{eqnarray}
and
\begin{eqnarray}
\label{det_err22}
I_{12}'\leq C\Delta t\sum_{i=0}^{m-1}\|X^h(t_i)-Y^h_i\|_{L^2(\Omega,\mathcal{H})}.
\end{eqnarray}
Adding \eqref{det_err21} and \eqref{det_err22}, we obtain 
\begin{eqnarray}
\label{det_err2}
I'_1\leq C\Delta t^{\frac{2H+\beta-1}2}+C\Delta t\sum_{i=0}^{m-1}\|X^h(t_i)-Y^h_i\|_{L^2(\Omega,\mathcal{H})}.
\end{eqnarray}
We estimate at now $I_2'$. Using triangle inequality and the estimate $(a+b)^2\leq 2a^2+2b^2$, we split it in three terms
\begin{eqnarray}
\label{stoerr'}
I_2'^2&=&\left\|\int_{0}^{t_m}S_h(t_m-s)P_h\phi(s)-S_h(t_m-[s])P_h\phi([s]) dB^H(s)\right\|^2_{L^2(\Omega,\mathcal{H})}\nonumber\\
&\leq& 2\left\|\int_{0}^{t_m}S_h(t_m-s)P_h\left(\phi(s)-\phi([s])\right)dB^H(s)\right\|^2_{L^2(\Omega,\mathcal{H})}\nonumber\\
&+&2\left\|\int_{0}^{t_m}\left(S_h(t_m-s)-S_h(t_m-[s])\right)P_h\phi([s]) dB^H(s)\right\|^2_{L^2(\Omega,\mathcal{H})}\nonumber\\
&\leq& 2\left\|\int_{0}^{t_m}S_h(t_m-s)P_h\left(\phi(s)-\phi([s])\right)dB^H(s)\right\|^2_{L^2(\Omega,\mathcal{H})}\nonumber\\
&+&2\left\|\int_{0}^{t_m}\left(S_h(t_m-s)-S_h(t_m-[s])\right)P_h\left(\phi([s])-\phi(t_{m-1})\right) dB^H(s)\right\|^2_{L^2(\Omega,\mathcal{H})}\nonumber\\
&+&2\left\|\int_{0}^{t_m}\left(S_h(t_m-s)-S_h(t_m-[s])\right)P_h\phi(t_{m-1}) dB^H(s)\right\|^2_{L^2(\Omega,\mathcal{H})}\nonumber\\&=:& 2I_{21}'^2+4I_{22}'^2+4I_{23}'^2.
\end{eqnarray}
Thanks to \eqref{stoerr1} we have
\begin{eqnarray}
\label{stoerr1'}
I_{21}'^2\leq C\Delta t^{2H+\beta-1}.
\end{eqnarray}
Thereafter, \eqref{bound}, the change of variable $j=m-k$ and $\varsigma=t_m-s$, inserting an appropriate power of $A_h$, \cite[(81)]{Mukb} and Remark \ref{disc_semigroup} (more precisely \eqref{semigroup_prp1} with $\gamma=\frac{2H+\beta-1}2$, \eqref{semigroup_prp2} and \assref{noise} (more precisely \eqref{ass_noise_term3}) yields 

\begin{eqnarray}
\label{stoerr2'}
I_{22}'^2&=&\left\|\int_{0}^{t_m}\left(S_h(t_m-s)-S_h(t_m-[s])\right)P_h\left(\phi([s])-\phi(t_{m-1})\right) dB^H(s)\right\|^2_{L^2(\Omega,\mathcal{H})}\nonumber\\
&\leq&C\int_{0}^{t_m}\left\|\left(S_h(t_m-s)-S_h(t_m-[s])\right)P_h\left(\phi([s])-\phi(t_{m-1})\right)\right\|^2_{L^0_2} ds\nonumber\\
&\leq&C\sum_{k=0}^{m-1}\int_{t_k}^{t_{k+1}}\left\|\left(S_h(t_m-s)-S_h(t_m-t_k)\right)P_h\left(\phi(t_k)-\phi(t_{m-1})\right)\right\|^2_{L^0_2} ds\nonumber\\
&\leq&C\sum_{j=1}^{m}\int_{t_{j-1}}^{t_j}\left\|\left(S_h(\varsigma)-S_h(t_j)\right)A_h^{\frac{1-\beta}2}\right\|^2_{L(\mathcal{H})}\nonumber\\
&&\times\left\|A_h^{\frac{\beta-1}2}P_h\left(\phi(t_{m-j})-\phi(t_{m-1})\right)\right\|^2_{L^0_2} d\varsigma\nonumber\\
&\leq&C\sum_{j=1}^{m}\int_{t_{j-1}}^{t_j}\left\|A_h^{-\frac{2H+\beta-1}2}\left(S_h(t_j-\varsigma)-I\right)\right\|^2_{L(\mathcal{H})}\left\|A_h^HS_h(\varsigma)\right\|^2_{L(\mathcal{H})}\nonumber\\
&&\times\left\|A^{\frac{\beta-1}2}\left(\phi(t_{m-j})-\phi(t_{m-1})\right)\right\|^2_{L^0_2} d\varsigma\nonumber\\
&\leq&C\sum_{j=1}^{m}\int_{t_{j-1}}^{t_j}\Delta t^{2H+\beta-1}\varsigma^{-2H}\left\|A^{\frac{\beta-1}2}\left(\phi(t_{m-j})-\phi(t_{m-1})\right)\right\|^2_{L^0_2} d\varsigma\nonumber\\
&\leq&C\Delta t^{2H+\beta-1}\left(\sum_{j=1}^{m}\int_{t_{j-1}}^{t_j}\varsigma^{-2H}t_{j-1}^{2\delta} d\varsigma\right)\nonumber\\
&\leq&C\Delta t^{2H+\beta-1}\left(\sum_{j=1}^{m}\int_{t_{j-1}}^{t_j}\varsigma^{2\delta-2H} d\varsigma\right)\nonumber\\
&\leq&C\Delta t^{2H+\beta-1}\left(\int_{0}^{t_m}\varsigma^{2\delta-2H} d\varsigma\right)\nonumber\\
&\leq&C\Delta t^{2H+\beta-1}t_m^{2\delta-2H+1}\leq C\Delta t^{2H+\beta-1}.
\end{eqnarray}
Finally, using \eqref{stoint1}, inserting an appropriate power of $A_h$, \cite[(81)]{Mukb} and Remark \ref{disc_semigroup} (more precisely \eqref{semigroup_prp1} with $\gamma=\frac{2H+\beta-1}2$, \eqref{semigroup_prp2} and \eqref{semigroup_prp5}), Assumption \ref{noise} ( more precisely \eqref{ass_noise_term2}) we obtain

\begin{eqnarray}
\label{stoerr3'}
I_{23}'^2&\leq&\left\|\int_{0}^{t_m}\left(S_h(t_m-s)-S_h(t_m-[s])\right)P_h\phi(t_{m-1}) dB^H(s)\right\|_{L^2(\Omega,\mathcal{H})}^2\nonumber\\
&\leq& C_H\sum_{i\in\N^d}\left(\int_{0}^{t_m}\left\|\left(S_h(t_m-s)-S_h(t_m-[s])\right)P_h\phi(t_{m-1}) Q^{\frac 12}e_i\right\|^{\frac 1H}ds\right)^{2H}\nonumber\\
&\leq& C_H\left(\int_{0}^{t_m}\left\|(S_h(t_m-s)-S_h(t_m-[s]))A_h^{\frac{1-\beta}2}\right\|^{\frac 1H}_{L(\mathcal{H})}ds\right)^{2H}\nonumber\\
&&\times\left(\sum_{i\in\N^d}\left\|A_h^{\frac{\beta-1}2}P_h\phi(t_{m-1}) Q^{\frac 12}e_i\right\|^2\right)\nonumber\\
&\leq& C_H\left(\int_{0}^{t_m}\left\|A_h^{-\frac{2H+\beta-1}2}(S_h(s-[s])-I)\right\|^{\frac 1H}_{L(\mathcal{H})}\left\|A_h^{H}S_h(t_m-s)\right\|^{\frac 1H}_{L(\mathcal{H})}ds\right)^{2H}\nonumber\\
&&\hspace{2cm}\left(\sum_{i\in\N^d}\left\|A^{\frac{\beta-1}2}\phi(t_{m-1}) Q^{\frac 12}e_i\right\|^2\right)\nonumber\\
&\leq& C_H\left(\int_{0}^{t_m}\left\|A_h^{H}S_h(t_m-s)\right\|^{\frac 1H}_{L(\mathcal{H})}(s-[s])^{\frac{2H+\beta-1}{2H}}ds\right)^{2H}\left\|A^{\frac{\beta-1}2}\phi(t_{m-1})\right\|^2_{L^0_2}\nonumber\\
&\leq& C_H\Delta t^{2H+\beta-1}\left(\int_{0}^{t_m}\left\|A_h^{H}S_h(t_m-s)\right\|^{\frac 1H}_{L(\mathcal{H})}ds\right)^{2H}\left\|A^{\frac{\beta-1}2}\phi(t_{m-1})\right\|^2_{L^0_2}\nonumber\\
&\leq& C\Delta t^{2H+\beta-1}\left\|A^{\frac{\beta-1}2}\phi(t_{m-1})\right\|^2_{L^0_2}\nonumber\\
&\leq& C\Delta t^{2H+\beta-1},
\end{eqnarray}
hence inserting \eqref{stoerr1'}-\eqref{stoerr3'} in \eqref{stoerr'} and taking the square-root gives 
\begin{equation}
\label{sto_err2'}
I_2'\leq C\Delta t^{\frac{2H+\beta-1}2}.
\end{equation}
Adding \eqref{det_err2} and \eqref{sto_err2'} we obtain
\begin{eqnarray}
\label{err2'}
err_2\leq C\Delta t^{\frac{2H+\beta-1}2}+C\Delta t\sum_{i=0}^{m-1}\|X^h(t_i)-Y^h_i\|_{L^2(\Omega,\mathcal{H})}.
\end{eqnarray}
Using the discrete version of the Gronwall inequality yields
\begin{eqnarray}
\label{time_err2}
err_2\leq C\Delta t^{\frac{2H+\beta-1}2}.
\end{eqnarray}
Combining \eqref{spa_err} and \eqref{time_err2} completes the proof.$\hfill\square$
\subsection{Proof of Theorem \ref{cstrongconvthm} for SERS scheme}
Before moving  to the proof, we first present some preparatory results. Thanks to \assref{der} and the works done in \cite{Muka} we obtain
\begin{lem}\cite[Lemma 5]{Muka}
\label{lemSERS1}
For all $M\in\N$ and all $\omega\in\Omega$, there is a positive constant $C'_1$ independent of $h$, $m$, $\Delta t$ and the sample $\omega$ such that
\begin{eqnarray*}
\left\|e^{(-A_h+J^h_m(\omega))t}\right\|_{L(\mathcal{H})}\leq C'_1,\hspace{2cm}0\leq t\leq T.
\end{eqnarray*}
\end{lem}
\begin{lem}\cite[Lemma 6]{Muka}
\label{lemSERS2}
The function $G^h_m(\omega)$ defined by \eqref{remain} satisfies the global Lipschitz condition 
with a uniform constant $C'>0$, independent of $h$, $m$ and $\omega$ such that
\begin{eqnarray*}
\left\|G^h_m(\omega)(u^h)-G^h_m(\omega)(v^h)\right\|\leq C'\|u^h-v^h\|\hspace{2cm}\forall u^h,\hspace{0.05cm}v^h\in V_h.
\end{eqnarray*}
\end{lem}
\begin{lem}\cite[Lemma 9]{Muka}
\label{lemSERS3}
For all $\omega\in\Omega$, the stochastic perturbed semigroup $S^h_m(\omega)(t):=e^{(-A_h+J^h_m(w))t}$ satisfies the following properties
\begin{enumerate}
\item [(i)] For $\gamma_1,\gamma_2\leq 1$ such that $0\leq\gamma_1+\gamma_2\leq 1$, 
\begin{eqnarray*}
\left\|A_h^{-\gamma_1}(S^h_m(\omega)(t)-I)A_h^{-\gamma_2}\right\|_{L(\mathcal{H})}\leq C t^{\gamma_1+\gamma_2},\hspace{2cm}t\in(0,T].
\end{eqnarray*}
\item[(ii)] For $\gamma_1\geq 0$ we have
\begin{eqnarray*}
\left\|S^h_m(\omega)(t)A_h^{\gamma_1}\right\|_{L(\mathcal{H})}\leq C t^{-\gamma_1},\hspace{2cm}t\in(0,T].
\end{eqnarray*}
\item[(iii)] For $\gamma_1\geq 0$ and $0\leq\gamma_2<1$such that $\gamma_2-\gamma_1\geq 0$, we have 
\begin{eqnarray*}
\|A_h^{-\gamma_1}S^h_m(\omega)(t)A_h^{\gamma_2}\|_{L(\mathcal{H})}\leq C t^{\gamma_1-\gamma_2},\hspace{2cm}t\in(0,T].
\end{eqnarray*}
\item[(iv)] For $\gamma_1,\gamma_2>0$ such that $0\leq\gamma_1-\gamma_2\leq 1$, then the following estimate holds 
\begin{eqnarray*}
\|A_h^{-\gamma_1}(S^h_m(\omega)(t)-I)A_h^{\gamma_2}\|_{L(\mathcal{H})}\leq C t^{\gamma_1-\gamma_2},\hspace{2cm}t\in(0,T].
\end{eqnarray*}
\end{enumerate} 
\end{lem}
\begin{lem}\cite[Lemma 10]{Muka}
\label{lemSERS4}
The stochastic perturbed semigroup $S^h_m(\omega)$ satisfies the following property
\begin{eqnarray*}
\left\|e^{(-A_h+J^h_m(\omega))\Delta t}\cdot\cdot\cdot e^{(-A_h+J^h_k(\omega))\Delta t}A_h^{\nu}\right\|_{L(\mathcal{H})}\leq C t_{m+1-k}^{-\nu},\quad 0\leq \nu<1.
\end{eqnarray*}
Where $C$ is a positive constant independent of $m$, $k$, $h$, $\Delta t$ and the sample $\omega$.
\end{lem}

We can now prove our theorem. As in the proof of the previous schemes, we split the fully discrete error in two terms as
\begin{eqnarray}
\label{errSERS}
\|X(t_m)-Z^h_m\|_{L^2(\Omega,\mathcal{H})}&\leq& \|X(t_m)-X^h(t_m)\|_{L^2(\Omega,\mathcal{H})}+\|X^h(t_m)-Z^h_m\|_{L^2(\Omega,\mathcal{H})}\nonumber\\
&=:& err_0+err_3.
\end{eqnarray}
By \lemref{serror} we have the estimate of  the space error
\begin{equation*}
err_0\leq C h^{2H+\beta-1}.
\end{equation*}
Consider now the estimate of the time error $err_3$. We recall that the semidiscrete problem \eqref{pb2} can be rewritten as
\begin{equation}
\label{pb3}
dX^h(t)=[-A_hX^h(t)+J^h_mX^h(t)+G^h_m(X^h(t))]dt+P_h\phi(t) dB^H(t)
\end{equation}
for all $t_m\leq t\leq t_{m+1}$ where $J^h_m$ and $G^h_m$ is given by \eqref{fre_der} and \eqref{remain}. Hence the exact solution at time $t_m$ of the semidiscrete problem \eqref{pb3} is given by
\begin{eqnarray}
\label{mild_sol_pb3}
X^h(t_m)&=&e^{(-A_h+J^h_{m-1})\Delta t}X^h(t_{m-1})+\int_{t_{m-1}}^{t_m}e^{(-A_h+J^h_{m-1})(t_m-s)}G^h_{m-1}(X^h(s))ds\nonumber\\
&&+\int_{t_{m-1}}^{t_m}e^{(-A_h+J^h_{m-1})(t_m-s)}P_h\phi(s) dB^H(s)
\end{eqnarray}
and the numerical solution \eqref{SERS} can be rewritten as 
\begin{eqnarray}
\label{SERS1}
Z^h_m&=&e^{(-A_h+J^h_{m-1})\Delta t}Z^h_{m-1}+\int_{t_{m-1}}^{t_m}e^{(-A_h+J^h_{m-1})(t_m-s)}G^h_{m-1}(Z^h_{m-1})ds\nonumber\\
&&+\int_{t_{m-1}}^{t_m}e^{(-A_h+J^h_{m-1})\Delta t}P_h\phi(t_{m-1}) dB^H(s).
\end{eqnarray}
If $m=1$ then from \eqref{mild_sol_pb3} and \eqref{SERS1} we obtain 
\begin{eqnarray}
\label{SERSerrm=1}
&&\|X^h(t_1)-Z^h_1\|_{L^2(\Omega,\mathcal{H})}\nonumber\\
&\leq& \left\|\int_{0}^{\Delta t}e^{(-A_h+J^h_0)(\Delta t-s)}(G^h_{0}(X^h(s))-G^h_{0}(Z^h_0))ds\right\|_{L^2(\Omega,\mathcal{H})}\nonumber\\
&&+\left\|\int_{0}^{\Delta t}e^{(-A_h+J^h_0)(\Delta t-s)}P_h\phi(s)-e^{(-A_h+J^h_0)\Delta t}P_h\phi(0) dB^H(s)\right\|_{L^2(\Omega,\mathcal{H})}\nonumber\\
&=:&I''+II''.
\end{eqnarray}
By \cite[(93)]{Muka} we have the estimate
\begin{equation}
\label{det_err3_m=1}
I''\leq C\Delta t.
\end{equation}
Using triangle inequality and the estimate $(a+b)^2\leq 2a^2+2b^2$ we split it in three terms
\begin{eqnarray}
\label{stoerr"m=1}
I_2'^2&=&\left\|\int_{0}^{t_m}e^{(-A_h+J^h_0)(\Delta t-s)}P_h\phi(s)-e^{(-A_h+J^h_0)\Delta t}P_h\phi(0) dB^H(s)\right\|^2_{L^2(\Omega,\mathcal{H})}\nonumber\\
&\leq& 2\left\|\int_{0}^{t_m}e^{(-A_h+J^h_0)(\Delta t-s)}P_h\left(\phi(s)-\phi(0)\right)dB^H(s)\right\|^2_{L^2(\Omega,\mathcal{H})}\nonumber\\
&&+2\left\|\int_{0}^{t_m}\left(e^{(-A_h+J^h_0)(\Delta t-s)}-e^{(-A_h+J^h_0)\Delta t}\right)P_h\phi(0) dB^H(s)\right\|^2_{L^2(\Omega,\mathcal{H})}\nonumber\\
&=:& 2II_1''^2+2II_2''^2.
\end{eqnarray}
For the estimate $II_1''^2$, using \eqref{bound}, inserting an appropriate power of $A_h$, \cite[(81)]{Mukb}, \assref{noise} (more precisely \eqref{ass_noise_term3}) and \lemref{lemSERS3} $(ii)$ with $\gamma_1=\frac{1-\beta}2$ we obtain
\begin{eqnarray}
\label{stoerr"m=1_1}
II_1''^2&=&\left\|\int_{0}^{\Delta t}e^{(-A_h+J^h_0)(\Delta t-s)}P_h\left(\phi(s)-\phi(0)\right)dB^H(s)\right\|^2_{L^2(\Omega,\mathcal{H})}\nonumber\\
&\leq& C\int_{0}^{\Delta t}\left\|e^{(-A_h+J^h_0)(\Delta t-s)}P_h\left(\phi(s)-\phi(0)\right)\right\|^2_{L^0_2}ds\nonumber\\
&\leq& C\int_{0}^{\Delta t}\left\|e^{(-A_h+J^h_0)(\Delta t-s)}A_h^{\frac{1-\beta}2}\right\|^2_{L(\mathcal{H})}\left\|A_h^{\frac{\beta-1}2}P_h\left(\phi(s)-\phi(0)\right)\right\|^2_{L^0_2}ds\nonumber\\
&\leq& C\int_{0}^{\Delta t}(\Delta t-s)^{\beta-1}s^{2\delta}ds\nonumber\\
&\leq& C\Delta t^{2H+\beta-1}\left(\int_{0}^{\Delta t}(\Delta t-s)^{\beta-1}ds\right)\nonumber\\
&\leq& C\Delta t^{2H+\beta-1}.
\end{eqnarray}
We denote by $\epsilon$ a positive constant small enough, using \eqref{stoint1}, inserting an appropriate power of $A_h$, Assumption \ref{noise}, \cite[(83)]{Mukb}, \lemref{lemSERS3} $(ii)$ with $\gamma_1=H-\frac{\epsilon}2$, $(iv)$ with $\gamma_1=H-\frac{\epsilon}2$ and $\gamma_2=\frac{1-\beta}2$ if $0\leq \beta < 1$ ( or $(i)$ with $\gamma_1=H-\frac{\epsilon}2$ and $\gamma_2=0$ if $\beta=1$) we have
{\small
\begin{eqnarray}
\label{stoerr"m=1_2}
II_2''^2&=&\mathbb{E}\left\|\int_{0}^{\Delta t}\left(e^{(-A_h+J^h_0)(\Delta t-s)}-e^{(-A_h+J^h_0)\Delta t}\right)P_h\phi(0) dB^H(s)\right\|^2\nonumber\\
&\leq& C_H\sum_{i\in\N^d}\left(\int_{0}^{\Delta t}\left\|(e^{(-A_h+J^h_0)(\Delta t-s)}-e^{(-A_h+J^h_0)\Delta t})P_h\phi(0) Q^{\frac 12}e_i\right\|^{\frac 1H}ds\right)^{2H}\nonumber\\
&\leq& C_H\sum_{i\in\N^d}\left(\int_{0}^{\Delta t}\left\|e^{(-A_h+J^h_0)(\Delta t-s)}\left(e^{(-A_h+J^h_0)s}-I\right)A_h^{\frac{1-\beta}2}\right\|^{\frac 1H}_{L(\mathcal{H})}\right.\nonumber\\
&&\left.\left\|A_h^{\frac{\beta-1}2}P_h\phi(0) Q^{\frac 12}e_i\right\|^{\frac 1H}ds\right)^{2H}\nonumber\\
&\leq& C_H\left(\int_{0}^{\Delta t}\left\|e^{(-A_h+J^h_0)(\Delta t-s)}A_h^{H-\frac{\epsilon}2}\right\|^{\frac 1H}_{L(\mathcal{H})}\left\|A_h^{\frac{\epsilon}2-H}\left(e^{(-A_h+J^h_0)s}-I\right)A_h^{\frac{1-\beta}2}\right\|^{\frac 1H}_{L(\mathcal{H})}ds\right)^{2H}\nonumber\\
&&\hspace{2cm}\left(\sum_{i\in\N^d}\left\|A^{\frac{\beta-1}2}\phi(0) Q^{\frac 12}e_i\right\|^2\right)\nonumber\\
&\leq& C_H\left(\int_{0}^{\Delta t}(\Delta t-s)^{\frac{\epsilon}{2H}-1}s^{\frac{2H+\beta-1-\epsilon}{2H}}ds\right)^{2H}\left\|A^{\frac{\beta-1}2}\phi(0)\right\|^2_{L^0_2}\nonumber\\
&\leq& C\Delta t^{2H+\beta-1-\varepsilon}\left(\int_{0}^{\Delta t}(\Delta t-s)^{\frac{\epsilon}{2H}-1}ds\right)^{2H}\nonumber\\
&\leq& C\Delta t^{2H+\beta-1-\epsilon}\Delta t^{\epsilon}\nonumber\\
&\leq& C\Delta t^{2H+\beta-1}.
\end{eqnarray}
}
Hence putting \eqref{stoerr"m=1_2} and \eqref{stoerr"m=1_2} in \eqref{stoerr"m=1} and taking the square-root gives 
\begin{eqnarray}
\label{sto_err3_m=1}
II''\leq C\Delta t^\frac{\beta+2H-1}2.
\end{eqnarray}
Adding \eqref{det_err3_m=1} and \eqref{sto_err3_m=1} yields
\begin{eqnarray}
\label{time_err3_m=1}
\|X^h(t_1)-Z^h_1\|_{L^2(\Omega,\mathcal{H})}\leq C\Delta t^\frac{\beta+2H-1}2.
\end{eqnarray}

For $m\geq 2$, we recall that the solution at $t_m$ of the semidiscrete problem \eqref{pb3} is given by
{\small
\begin{eqnarray}
\label{mild_sol_pb3'}
&&X^h(t_m)\\
&=&e^{(-A_h+J^h_{m-1})\Delta t}\cdot\cdot\cdot e^{(-A_h+J^h_{m-1})\Delta t}X^h(0)\nonumber\\
&+&\int_{t_{m-1}}^{t_m}e^{(-A_h+J^h_{m-1})(t_m-s)}G^h_{m-1}(X^h(s))ds\nonumber\\
&+&\int_{t_{m-1}}^{t_m}e^{(-A_h+J^h_{m-1})(t_m-s)}P_h\phi(s) dB^H(s)\nonumber\\
&+&\int_{0}^{t_{m-1}}e^{(-A_h+J^h_{m-1})\Delta t}\cdot\cdot\cdot e^{(-A_h+J^h_{[s]^m+1})\Delta t}e^{(-A_h+J^h_{[s]^m})([s]+\Delta t-s)}G^h_{[s]^m}(X^h(s))ds\nonumber\\
&+&\int^{t_{m-1}}_0e^{(-A_h+J^h_{m-1})\Delta t}\cdot\cdot\cdot e^{(-A_h+J^h_{[s]^m+1})\Delta t}e^{(-A_h+J^h_{[s]^m})([s]+\Delta t-s)}P_h\phi(s) dB^H(s).\nonumber
\end{eqnarray} 
}
We recall also that the numerical solution at $t_m$ given by \eqref{SERS1} can be rewritten as
{\small
\begin{eqnarray}
\label{SERS1'}
Z^h_m&=&e^{(-A_h+J^h_{m-1})\Delta t}\cdot\cdot\cdot e^{(-A_h+J^h_{m-1})\Delta t}X^h(0)\nonumber\\
&+&\int_{t_{m-1}}^{t_m}e^{(-A_h+J^h_{m-1})(t_m-s)}G^h_{m-1}(Z^h_{m-1})ds\nonumber\\
&+&\int_{t_{m-1}}^{t_m}e^{(-A_h+J^h_{m-1})\Delta t}P_h\phi(t_{m-1}) dB^H(s)\nonumber\\
&+&\int_{0}^{t_{m-1}}e^{(-A_h+J^h_{m-1})\Delta t}\cdot\cdot\cdot e^{(-A_h+J^h_{[s]^m+1})\Delta t}e^{(-A_h+J^h_{[s]^m})([s]+\Delta t-s)}G^h_{[s]^m}(Z^h_{[s]^m})ds\nonumber\\
&+&\int_{0}^{t_{m-1}}e^{(-A_h+J^h_{m-1})\Delta t}\cdot\cdot\cdot e^{(-A_h+J^h_{[s]^m+1})\Delta t}e^{(-A_h+J^h_{[s]^m})\Delta t}P_h\phi([s]) dB^H(s).
\end{eqnarray} 
}
Using \eqref{mild_sol_pb3'}, \eqref{SERS1'} and the triangle inequality, we have
\begin{eqnarray}
\label{SERSmneq1}
\|X^h(t_m)-Z^h_m\|_{L^2(\Omega,\mathcal{H})}\leq III''+IV''+V''+VI''
\end{eqnarray}
where 
{\small
\begin{eqnarray*}
&&III''=\left\|\int_{t_{m-1}}^{t_m}e^{(-A_h+J^h_{m-1})(t_m-s)}[G^h_{m-1}(X^h(s))-G^h_{m-1}(Z^h_{m-1})]ds\right\|_{L^2(\Omega,\mathcal{H})},\\
&&IV''=\left\|\int_{t_{m-1}}^{t_m}e^{(-A_h+J^h_{m-1})(t_m-s)}P_h\phi(s)-e^{(-A_h+J^h_{m-1})\Delta t}P_h\phi(t_{m-1}) dB^H(s)\right\|_{L^2(\Omega,\mathcal{H})},
\end{eqnarray*}
}
\begin{eqnarray*}
V''&=&\left\|\int_0^{t_{m-1}}e^{(-A_h+J^h_{m-1})\Delta t}\cdot\cdot\cdot e^{(-A_h+J^h_{[s]^m+1})\Delta t}e^{(-A_h+J^h_{[s]^m})([s]+\Delta t-s)}\right.\nonumber\\
&&\left.\left[G^h_{[s]^m}(X^h(s))-G^h_{[s]^m}(Z^h_{[s]^m})\right]ds\right\|_{L^2(\Omega,\mathcal{H})},
\end{eqnarray*}
and 
\begin{eqnarray*}
VI''&=&\left\|\int_0^{t_{m-1}}e^{(-A_h+J^h_{m-1})\Delta t}\cdot\cdot\cdot e^{(-A_h+J^h_{[s]^m+1})\Delta t}e^{(-A_h+J^h_{[s]^m})([s]+\Delta t-s)}P_h\phi(s)\right.\\
&-&\left.e^{(-A_h+J^h_{m-1})\Delta t}\cdot\cdot\cdot e^{(-A_h+J^h_{[s]^m+1})\Delta t}e^{(-A_h+J^h_{[s]^m})\Delta t}P_h\phi([s]) dB^H(s)\right\|_{L^2(\Omega,\mathcal{H})}.
\end{eqnarray*}
By Lemma \ref{lemSERS1}, triangle inequality and Lemma \ref{lem3} \eqref{reg2}, we  easily have
\begin{eqnarray}
\label{det_err3'}
III''\leq C\Delta t^{\frac{2H+\beta-1}2}+C\Delta t\|X^h(t_{m-1})-Z^h_{m-1}\|_{L^2(\Omega,\mathcal{H})}.
\end{eqnarray}
In a similar way, using Lemma \ref{lemSERS4} with $\nu=0$, \lemref{lemSERS1}, \lemref{lemSERS2}, triangle inequality and \lemref{lem3} \eqref{reg2}, we obtain
\begin{eqnarray}
\label{det_err3''}
V''\leq C\Delta t^{\frac{2H+\beta-1}2}+C\Delta t\sum_{k=0}^{m-2}\|X^h(t_k)-Z^h_k\|_{L^2(\Omega,\mathcal{H})}.
\end{eqnarray}
For the estimate $IV''$, Using triangle inequality and the estimate $(a+b)^2\leq 2a^2+2b^2$ we split it in two terms
{\small
\begin{eqnarray}
\label{stoerr3m>1}
IV''^2&=&\left\|\int_{t_{m-1}}^{t_m}e^{(-A_h+J^h_{m-1})(t_m-s)}P_h\phi(s)-e^{(-A_h+J^h_{m-1})\Delta t}P_h\phi(t_{m-1}) dB^H(s)\right\|^2_{L^2(\Omega,\mathcal{H})}\nonumber\\
&\leq& 2\left\|\int_{0}^{t_m}e^{(-A_h+J^h_{m-1})(t_m-s)}P_h\left(\phi(s)-\phi(t_{m-1})\right)dB^H(s)\right\|^2_{L^2(\Omega,\mathcal{H})}\nonumber\\
&&+2\left\|\int_{0}^{t_m}\left[e^{(-A_h+J^h_{m-1})(t_m-s)}-e^{(-A_h+J^h_{m-1})\Delta t}\right]P_h\phi(t_{m-1}) dB^H(s)\right\|^2_{L^2(\Omega,\mathcal{H})}\nonumber\\
&=:& 2IV_1''^2+2IV_2''^2.
\end{eqnarray}
}
Firstly, using \eqref{bound}, inserting an appropriate power of $A_h$, \cite[(81)]{Mukb}, \assref{noise} (more precisely \eqref{ass_noise_term3}) and \lemref{lemSERS3} $(ii)$ with $\gamma_1=\frac{1-\beta}2$ we obtain
\begin{eqnarray}
\label{stoerr3m>1_1}
IV_1''^2&=&\left\|\int_{0}^{\Delta t}e^{(-A_h+J^h_{m-1})(t_m-s)}P_h\left(\phi(s)-\phi(t_{m-1})\right)dB^H(s)\right\|^2_{L^2(\Omega,\mathcal{H})}\nonumber\\
&\leq& C\int_{0}^{\Delta t}\left\|e^{(-A_h+J^h_{m-1})(t_m-s)}P_h\left(\phi(s)-\phi(t_{m-1})\right)\right\|^2_{L^0_2}ds\nonumber\\
&\leq& C\int_{0}^{\Delta t}\left\|e^{(-A_h+J^h_{m-1})(t_m-s)}A_h^{\frac{1-\beta}2}\right\|^2_{L(\mathcal{H})}\left\|A_h^{\frac{\beta-1}2}P_h\left(\phi(s)-\phi(t_{m-1})\right)\right\|^2_{L^0_2}ds\nonumber\\
&\leq& C\int_{0}^{\Delta t}(t_m-s)^{\beta-1}\left\|A^{\frac{\beta-1}2}\left(\phi(s)-\phi(t_{m-1})\right)\right\|^2_{L^0_2}ds\nonumber\\
&\leq& C\int_{0}^{\Delta t}(t_m-s)^{\beta-1}(s-t_{m-1})^{2\delta}ds\nonumber\\
&\leq& C\Delta t^{2H+\beta-1}\left(\int_{0}^{\Delta t}(t_m-s)^{\beta-1}ds\right)\nonumber\\
&\leq& C\Delta t^{2H+\beta-1}.
\end{eqnarray}
Secondly by \eqref{stoint1}, inserting an appropriate power of $A_h$, Assumption \ref{noise}, \cite[(81)]{Mukb}, \lemref{lemSERS3} $(ii)$ with $\gamma_1=H-\frac{\epsilon}2$, $(iv)$ with $\gamma_1=H-\frac{\epsilon}2$ and $\gamma_2=\frac{1-\beta}2$, we obtain
{\small
\begin{eqnarray}
\label{stoerr3m>1_2}
IV^{''2}&=&\mathbb{E}\left\|\int_{t_{m-1}}^{t_m}\left(e^{(-A_h+J^h_{m-1})(t_m-s)}-e^{(-A_h+J^h_{m-1})\Delta t}\right)P_h\phi(t_{m-1}) dB^H(s)\right\|^2\nonumber\\
&\leq& C_H\sum_{i\in\N^d}\left(\int_{t_{m-1}}^{t_m}\left\|\left[e^{(-A_h+J^h_{m-1})(t_m-s)}-e^{(-A_h+J^h_{m-1})\Delta t}\right]P_h\phi(t_{m-1}) Q^{\frac 12}e_i\right\|^{\frac 1H}ds\right)^{2H}\nonumber\\
&\leq& C_H\left(\int_{t_{m-1}}^{t_m}\left\|e^{(-A_h+J^h_{m-1})(t_m-s)}\left(e^{(-A_h+J^h_{m-1})(s-t_{m-1})}-I\right)A_h^{\frac{1-\beta}2}\right\|^{\frac 1H}_{L(\mathcal{H})}ds\right)^{2H}\nonumber\\
&&\hspace{2cm}\left(\sum_{i\in\N^d}\left\|A^{\frac{\beta-1}2}\phi(t_{m-1}) Q^{\frac 12}e_i\right\|^2\right)\nonumber\\
&\leq& C_H\left(\int_{t_{m-1}}^{t_m}\left\|e^{(-A_h+J^h_{m-1})(t_m-s)}A_h^{H-\frac{\epsilon}2}\right\|^{\frac 1H}_{L(\mathcal{H})}\right.\nonumber\\
&&\left.\left\|A_h^{-H+\frac{\epsilon}2}\left(e^{(-A_h+J^h_{m-1})(s-t_{m-1})}-I\right)A_h^{\frac{1-\beta}2}\right\|^{\frac 1H}_{L(\mathcal{H})}ds\right)^{2H}\nonumber\\
&&\hspace{2cm}\left(\sum_{i\in\N^d}\left\|A^{\frac{\beta-1}2}\phi(t_{m-1}) Q^{\frac 12}e_i\right\|^2\right)\nonumber\\
&\leq& C_H\left(\int_{t_{m-1}}^{t_m}(t_m-s)^{-1+\frac{\epsilon}{2H}}(s-t_{m-1})^{\frac{2H+\beta-1-\epsilon}{2H}}ds\right)^{2H}\left\|A^{\frac{\beta-1}2}\phi(t_{m-1})\right\|^2_{L^0_2}\nonumber\\
&\leq& C\Delta t^{2H+\beta-1-\epsilon}\left(\int_{t_{m-1}}^{t_m}(t_m-s)^{-1+\frac{\epsilon}{2H}}ds\right)^{2H}\nonumber\\
&\leq& C\Delta t^{2H+\beta-1-\epsilon}\Delta t^{\epsilon}\nonumber\\
&\leq&  C\Delta t^{2H+\beta-1}.
\end{eqnarray}
}
Hence substituting \eqref{stoerr3m>1_1} and \eqref{stoerr3m>1_2} in \eqref{stoerr3m>1} and taking the square-root gives 
\begin{eqnarray}
\label{sto_err3'}
IV''\leq C\Delta t^{\frac{2H+\beta-1}2}.
\end{eqnarray}
For estimate $VI''$, using triangle inequality and the estimate $(a+b+c)^2\leq 3a^2+3b^2+3c^2$ we split it in two terms as

{\small
\begin{eqnarray}
\label{stoerr4m>1}
VI''^2
&\leq&3\left\|\int_0^{t_{m-1}}e^{(-A_h+J^h_{m-1})\Delta t}\cdot\cdot\cdot e^{(-A_h+J^h_{[s]^m+1})\Delta t}e^{(-A_h+J^h_{[s]^m})([s]+\Delta t-s)}\right.\nonumber\\
&&\left.P_h(\phi(s)-\phi([s]))dB^H(s)\right\|^2_{L^2(\Omega,\mathcal{H})}\nonumber\\
&+&3\|\int_0^{t_{m-1}}e^{(-A_h+J^h_{m-1})\Delta t}\cdot\cdot\cdot e^{(-A_h+J^h_{[s]^m+1})\Delta t}\nonumber\\
&&\left(e^{(-A_h+J^h_{[s]^m})([s]+\Delta t-s)}-e^{(-A_h+J^h_{[s]^m})\Delta t}\right)P_h\left(\phi([s])-\phi(t_{m-1})\right)dB^H(s)\|^2_{L^2(\Omega,\mathcal{H})}\nonumber\\
&+&3\|\int_0^{t_{m-1}}e^{(-A_h+J^h_{m-1})\Delta t}\cdot\cdot\cdot e^{(-A_h+J^h_{[s]^m+1})\Delta t}\nonumber\\
&&\left(e^{(-A_h+J^h_{[s]^m})([s]+\Delta t-s)}-e^{(-A_h+J^h_{[s]^m})\Delta t}\right)\nonumber\\
&&\times P_h\phi(t_{m-1})dB^H(s)\|^2_{L^2(\Omega,\mathcal{H})}\nonumber\\
&=:& 3VI_1''^2+3VI_2''^2+3VI_3''^2.
\end{eqnarray}
}
Let $\epsilon$ be a sufficient small number. At first, using \eqref{bound}, inserting an appropriate power of $A_h$, \cite[(81)]{Mukb}, \assref{noise} (more precisely \eqref{ass_noise_term3}), \lemref{lemSERS4} $(ii)$ with $\nu=H-\frac{\epsilon}2$, \lemref{lemSERS3} $(iii)$ with $\gamma_1=\gamma_2=\frac{1-\beta}{2}$, the variable change $j=m-k-1$ and \cite[(169)]{Muka} we have

\begin{eqnarray}
\label{stoerr4m>1_1}
VI_1''^2
&=&\left\|\int_0^{t_{m-1}}e^{(-A_h+J^h_{m-1})\Delta t}\cdot\cdot\cdot e^{(-A_h+J^h_{[s]^m+1})\Delta t}e^{(-A_h+J^h_{[s]^m})([s]+\Delta t-s)}\right.\nonumber\\
&&\left.P_h(\phi(s)-\phi([s]))dB^H(s)\right\|^2_{L^2(\Omega,\mathcal{H})}\nonumber\\
&\leq& C\int_{0}^{t_{m-1}}\left\|e^{(-A_h+J^h_{m-1})\Delta t}\cdot\cdot\cdot e^{(-A_h+J^h_{[s]^m+1})\Delta t}e^{(-A_h+J^h_{[s]^m})([s]+\Delta t-s)}\right.\nonumber\\
&&\left.P_h(\phi(s)-\phi([s]))\right\|_{L^0_2}^2ds\nonumber\\
&\leq& C\int_{0}^{t_{m-1}}\left\|e^{(-A_h+J^h_{m-1})\Delta t}\cdot\cdot\cdot e^{(-A_h+J^h_{[s]^m+1})\Delta t}A_h^{\frac{1-\beta}{2}}\right\|^2_{L(\mathcal{H})}\nonumber\\
&\times&\left\|A_h^{\frac{\beta-1}2}e^{(-A_h+J^h_{[s]^m})([s]+\Delta t-s)}A_h^{\frac{1-\beta}2}\right\|^2_{L(\mathcal{H})}\left\|A_h^{\frac{\beta-1}2}P_h\left(\phi(s)-\phi([s])\right)\right\|^2_{L^0_2}ds\nonumber\\
&\leq& C\int_{0}^{t_{m-1}}t_{m-[s]^m-1}^{-1+\beta}\left\|A^{\frac{\beta-1}2}\left(\phi(s)-\phi([s])\right)\right\|^2_{L^0_2}ds\nonumber\\
&\leq& C\int_{0}^{t_{m-1}}t_{m-[s]^m-1}^{-1+\beta}(s-[s])^{2\delta}ds\nonumber\\
&\leq& C\sum_{k=0}^{m-2}\int_{t_k}^{t_{k+1}}t_{m-k-1}^{-1+\beta}\Delta t^{2\delta}ds\nonumber\\
&\leq& C\Delta t^{2H+\beta-1}\left(\sum_{j=1}^{m-1}t_j^{-1+\beta}\Delta t\right)\nonumber\\
&\leq& C\Delta t^{2H+\beta-1}.
\end{eqnarray}
Using\eqref{bound}, inserting an appropriate power of $A_h$, \lemref{lemSERS4} with $\nu=H$, \lemref{lemSERS3} $(iii)$ with $\gamma_1=\gamma_2=H$, $(iv)$ with $\gamma_1=H$ and $\gamma_2=\frac {1-\beta}2$ if $0 < \beta < 1$ ( or $(i)$ with $\gamma_1=H$ and $\gamma_2=0$ if $\beta=1$), \cite[(81)]{Mukb}, \assref{noise} ( more precisely \eqref{ass_noise_term3}) the variable change $j=m-k-1$  and \cite[(169)]{Muka} yields

\begin{eqnarray}
\label{stoerr4m>1_2}
VI_2^{''2}
&=&\left\|\int_0^{t_{m-1}}e^{(-A_h+J^h_{m-1})\Delta t}\cdot\cdot\cdot e^{(-A_h+J^h_{[s]^m+1})\Delta t}\right.\nonumber\\
&&\left(e^{(-A_h+J^h_{[s]^m})([s]+\Delta t-s)}-e^{(-A_h+J^h_{[s]^m})\Delta t}\right)\nonumber\\
&&\hspace{4cm}\left.P_h\left(\phi([s])-\phi(t_{m-1})\right)dB^H(s)\right\|^2_{L^2(\Omega,\mathcal{H})}\nonumber\\
&\leq& C_H\sum_{i\in\N^d}\left(\int_0^{t_{m-1}}\left\|e^{(-A_h+J^h_{m-1})\Delta t}\cdot\cdot\cdot e^{(-A_h+J^h_{[s]^m+1})\Delta t}A_h^{H}\right\|^2_{L(\mathcal{H})}\right.\nonumber\\
&&\left\|A_h^{-H}e^{(-A_h+J^h_{[s]^m})([s]+\Delta t-s)}A_h^{H}\right\|^2_{L(\mathcal{H})}\nonumber\\
&&\left\|A_h^{-H}\left(e^{(-A_h+J^h_{[s]^m})(s-[s])}-I\right)A_h^{\frac{1-\beta}2}\right\|^2_{L(\mathcal{H})}\nonumber\\
&&\hspace{4cm}\left.\left\|A_h^{\frac{\beta-1}2}P_h\left(\phi([s])-\phi(t_{m-1})\right) Q^{\frac 12}e_i\right\|^2ds\right)\nonumber\\
&\leq& C_H\left(\int_0^{t_{m-1}}t^{-2H}_{m-[s]^m-1}(s-[s])^{2H+\beta-1}\left\|A^{\frac{\beta-1}2}\left(\phi(t_{m-1})-\phi([s])\right)\right\|^2_{L^0_2}ds\right)\nonumber\\
&\leq& C\Delta t^{2H+\beta-1}\left(\int_0^{t_{m-1}}t^{-2H}_{m-[s]^m-1}(t_{m-1}-[s])^{2\delta}ds\right)\nonumber\\
&\leq& C\Delta t^{2H+\beta-1}\left(\sum_{k=0}^{m-2}\int_{t_k}^{t_{k+1}}t^{-2H}_{m-k-1}(t_{m-1}-t_k)^{2H+\beta-1}ds\right)\nonumber\\
&\leq& C\Delta t^{2H+\beta-1}\left(\sum_{k=0}^{m-2}t^{-1+\beta}_{m-k-1}\Delta t\right)\nonumber\\
&\leq& C\Delta t^{2H+\beta-1}\left(\sum_{j=1}^{m-1}t^{-1+\beta}_{j}\Delta t\right)\nonumber\\
&\leq& C\Delta t^{2H+\beta-1}.
\end{eqnarray}
Afterwards, \eqref{stoint1}, inserting an appropriate power of $A_h$, \lemref{lemSERS4} with $\nu=H-\frac{\epsilon}2$, \lemref{lemSERS3} $(iii)$ with $\gamma_1=\gamma_2=H-\frac{\epsilon}2$, $(iv)$ with $\gamma_1=H-\frac{\epsilon}2$ and $\gamma_2=\frac {1-\beta}2$ if $0 < \beta < 1$ ( or $(i)$ with $\gamma_1=H-\frac{\epsilon}2$ and $\gamma_2=0$ if $\beta=1$), \cite[(81)]{Mukb}, \assref{noise} (more precisely \eqref{ass_noise_term2} ) the variable change $j=m-k-1$  and \cite[(169)]{Muka} yields

{\small
\begin{eqnarray}
\label{stoerr4m>1_3}
VI_3^{''2}&=&\left\|\int_0^{t_{m-1}}e^{(-A_h+J^h_{m-1})\Delta t}\cdot\cdot\cdot e^{(-A_h+J^h_{[s]^m+1})\Delta t}\right.\nonumber\\
&&\left(e^{(-A_h+J^h_{[s]^m})([s]+\Delta t-s)}-e^{(-A_h+J^h_{[s]^m})\Delta t}\right)P_h\phi(t_{m-1})dB^H(s)\|^2_{L^2(\Omega,\mathcal{H})}\nonumber\\
&\leq& C_H\sum_{i\in\N^d}\left(\int_0^{t_{m-1}}\left\|e^{(-A_h+J^h_{m-1})\Delta t}\cdot\cdot\cdot e^{(-A_h+J^h_{[s]^m+1})\Delta t}A_h^{H-\frac {\epsilon}2}\right\|^{\frac 1H}_{L(\mathcal{H})}\right.\nonumber\\
&&\left\|A_h^{\frac {\epsilon}2-H}e^{(-A_h+J^h_{[s]^m})([s]+\Delta t-s)}A_h^{H-\frac {\epsilon}2}\right\|^{\frac 1H}_{L(\mathcal{H})}\nonumber\\
&&\left.\left\|A_h^{\frac {\epsilon}2-H}\left(e^{(-A_h+J^h_{[s]^m})(s-[s])}-I\right)A_h^{\frac{1-\beta}2}\right\|^{\frac 1H}_{L(\mathcal{H})}\left\|A_h^{\frac{\beta-1}2}P_h\phi(t_{m-1}) Q^{\frac 12}e_i\right\|^{\frac 1H}ds\right)^{2H}\nonumber\\
&\leq& C_H\sum_{i\in\N^d}\left(\int_0^{t_{m-1}}t^{-1+\frac{\epsilon}{2H}}_{m-[s]^m-1}(s-[s])^{\frac{2H+\beta-1-\epsilon}{2H}}\left\|A^{\frac{\beta-1}2}\phi(t_{m-1}) Q^{\frac 12}e_i\right\|^{\frac 1H}ds\right)^{2H}\nonumber\\
&\leq& C\Delta t^{2H+\beta-1-\epsilon}\left(\int_0^{t_{m-1}}t^{-1+\frac{\epsilon}{2H}}_{m-[s]^m-1}ds\right)^{2H}\left(\sum_{i\in\N^d}\left\|A^{\frac{\beta-1}2}\phi(t_{m-1}) Q^{\frac 12}e_i\right\|^2\right)\nonumber\\
&\leq& C\Delta t^{2H+\beta-1-\epsilon}\left(\sum_{k=0}^{m-2}\int_{t_k}^{t_{k+1}}t^{-1+\frac{\epsilon}{2H}}_{m-k-1}ds\right)^{2H}\left\|A^{\frac{\beta-1}2}\phi(t_{m-1})\right\|^2_{L^0_2}\nonumber\\
&\leq& C\Delta t^{2H+\beta-1-\epsilon}\left(\sum_{k=0}^{m-2}t^{-1+\frac{\epsilon}{2H}}_{m-k-1}\Delta t\right)^{2H}\nonumber\\
&\leq& C\Delta t^{2H+\beta-1-\epsilon}\left(\sum_{j=1}^{m-1}t^{-1+\frac{\epsilon}{2H}}_{j}\Delta t\right)^{2H}\nonumber\\
&\leq& C\Delta t^{2H+\beta-1-\epsilon}.
\end{eqnarray}
}
Hence inserting \eqref{stoerr4m>1_1}, \eqref{stoerr4m>1_2} and \eqref{stoerr4m>1_3} in \eqref{stoerr4m>1} and taking the square-root gives 
\begin{eqnarray}
\label{sto_err3''}
VI''\leq C\Delta t^{\frac{2H+\beta-1-\epsilon}2}.
\end{eqnarray}
Adding \eqref{time_err3_m=1}, \eqref{det_err3'}, \eqref{det_err3''}, \eqref{sto_err3'}, \eqref{sto_err3''}
 and applying Gronwall's lemma yields
\begin{eqnarray}
\label{time_err3}
\|X^h(t_m)-Z^h_m\|_{L^2(\Omega,\mathcal{H})}\leq C\Delta t^{\frac{2H+\beta-1-\epsilon}2}.
\end{eqnarray}
Finally, adding \eqref{spa_err} and \eqref{time_err3} completes the proof.$\hfill\square$
\section{Extension  to SPDE driven simultaneously by fBm and Poisson random measure}
\label{convp}
\subsection{Numerical schemes}
Here  the goal is to show  how the previous results can be extended  to the following SPDE driven simultaneously by fBm and Poisson random measure.   The corresponding model equation is given by

{\small
\begin{eqnarray}
\label{model2}
dX(t, x)&=&\left[\nabla \cdot \left(\mathbf{D}\nabla X(t,x)\right)-\mathbf{q}\cdot\nabla X(t,x)+f(x,X(t,x))\right]dt+b(x, t))dB^H(t,x)\nonumber\\
&+&\int_{\chi}z_0\widetilde{N}(dz,dt),\quad x\in\Lambda,  \quad X(0)=X_0,\quad z_0\in\chi \quad t\in[0,T],
\end{eqnarray}
}
In the Hilbert space  $\mathcal{H}=L^2(\Lambda)$, \eqref{model2} is equivalent to \eqref{model3} where the linear operator $A$ and the nonlinear function $F$ are defined as in \eqref{operator} and \eqref{nemystskii}.
 The well posedness result  for $H=1/2$  presented  in  \cite{Alberverio} can easily  be extended  to \eqref{model3} for $H \in [1/2,1] $ by combining with \cite{Wanc}.
 The  corresponding  exponential Euler (SETD1)  scheme in integral  form is therefore given
\begin{eqnarray}
\label{PSTDE1}
Y^h_{m+1}&=&S_h(\Delta t)Y^h_m+\int_{t_m}^{t_{m+1}}S_h(t_{m+1}-s)P_hF(Y^h_m)ds\\
&+&\int_{t_m}^{t_{m+1}}S_h(\Delta t)P_h \phi(t_m) dB^H(s)+ \int_{t_m}^{t_{m+1}}\int_{\chi}S_h(\Delta t)P_h z_0\widetilde{N}(dz, ds),\nonumber
\end{eqnarray}
with $Y^h_0=P_hX_0$.
In the same way, the semi implicit scheme  is given by
\begin{eqnarray}
\label{IM1}
X^h_{m+1}&=&S_{h, \Delta t}Y^h_m+\int_{t_m}^{t_{m+1}}S_{h, \Delta t}P_hF(X^h_m)ds+\int_{t_m}^{t_{m+1}}S_{h, \Delta t}P_h \phi(t_m) dB^H(s)\nonumber\\
&+&\int_{t_m}^{t_{m+1}}\int_{\chi}S_{h, \Delta t}P_h z_0\widetilde{N}(dz, ds),
\end{eqnarray}
 while the  stochastic Exponential Rosenbrock Scheme (SERS) is given  by 
\begin{eqnarray}
\label{SERS_1}
Z^h_{m+1}&=&S_h^m(\Delta t)Z^h_m+\int_{t_m}^{t_{m+1}}S_h^m(t_{m+1}-s)P_hF(Z^h_m)ds\\
&+&\int_{t_m}^{t_{m+1}}S_h^m(\Delta t)P_h \phi(t_m) dB^H(s)+ \int_{t_m}^{t_{m+1}}\int_{\chi}S_h^m(\Delta t)P_h z_0\widetilde{N}(dz, ds),\nonumber
\end{eqnarray}
with $Z^h_0=P_hX_0$, 
where 
\begin{eqnarray}
\label{vend1a}
S_{h}^m(t):= e^{(-A_h+J_m^h)(t)}.
\end{eqnarray}
To obtain  the optimal order in time, as in \cite{manto}, we need the following assumption in Poisson measure noise.
\begin{ass}
 \label{assumption5}
 The covariance operator $Q:\mathcal{H}\longrightarrow \mathcal{H}$ and the jump coefficient  satisfy the following estimate
 \begin{eqnarray}
 \label{poisson}
 \Vert A^{\frac{\eta-1}{2}}z_0\Vert<\infty.
 \end{eqnarray}
 where $\eta=2H+\beta-1$ with $\beta\in(0,1]$ as in \assref{init} and \assref{noise} .
 \end{ass}
 \begin{remark}
 All the regularity results in space and time, both for continuous equation \eqref{model3} (or semi-discrete equation) important to achieve optimal convergence orders  can easily be extended  from  our results on  fBm  in \thmref{reg} and \lemref{lem3} by just following 
 \cite[Proposition 3.1]{manto}.
 \end{remark}
 
 \subsection{Convergence results for   SPDE with fBm and Poisson measure noise}
The convergence result is exactly as for fBm when \assref{assumption5}  is used
\begin{thm}
\label{cstrongconvthmp}
Let $X(t_m)$ be the mild solution of \eqref{model3} at time $t_m=m\Delta t$, $\Delta t\geq 0$.  Let $\zeta^h_m$ be the numerical approximations through \eqref{IM1} and \eqref{SERS1}($\zeta^h_m=X^h_m$ for implicit scheme, $\zeta^h_m=Z^h_m$ for SERS)
and $Y^h_m$ the numerical approximation through the  SETD1 given in \eqref{PSTDE1}. If Assumptions \ref{noise}-\ref{der} and \assref{assumption5} hold with $\beta\in(0,1]$, then 
\begin{equation}
\label{strongconvSETD1p}
\left(\mathbb{E}\|X(t_m)-Y^h_m\|^2\right)^{\frac 12}\leq C\left(h^{2H+\beta-1}+\Delta t^{\frac{2H+\beta-1}2}\right),
\end{equation}
and 
\begin{equation}
\label{strongconvimplandSERSp}
\left(\mathbb{E}\|X(t_m)-\zeta^h_m\|^2\right)^{\frac 12}\leq C\left(h^{2H+\beta-1}+\Delta t^{\frac{2H+\beta-1-\epsilon}2}\right),
\end{equation}
where $\epsilon$ is a positive constant small enough.
\end{thm}
\begin{cor}
Let $X(t_m)$ be the mild solution of \eqref{model3} (A self-adjoint) at time $t_m=m\Delta t$, $\Delta t\geq 0$.  Let $X^h_m$ be the numerical approximations through \eqref{IM1}. If Assumptions \ref{noise}-\ref{init} and \assref{assumption5} hold with $\beta\in(0,1]$, then 
\begin{equation}
\label{strongconvimplp}
\left(\mathbb{E}\|X(t_m)-X^h_m\|^2\right)^{\frac 12}\leq C\left(h^{2H+\beta-1}+\Delta t^{\frac{2H+\beta-1}2}\right).
\end{equation}
\end{cor}
\subsection{Proof of  convergence results  for  SPDE  with fBm and Poisson measure noise}
As in \cite{manto}, the proofs  are based on  Burkholder-Davis-Gundy Inequality where  the fBm version  
is given in \cite[Theorem 1.2 ]{G}. Under Assumptions \ref{noise}-\ref{init} and \assref{assumption5}, The regulatity result in time is
\begin{eqnarray}
\label{mild_sol5}
\|X^h(t_2)-X^h(t_1)\|_{L^2(\Omega,\mathcal{H})}\leq C(t_2-t_1)^{\frac{\min\left(2H+\beta-1,1\right)}2},\,0\leq t_1<t_2\leq T.
\end{eqnarray}
Where $C=C(\beta,L,T,H)$ is a positive constant and $\beta$ is the regularity parameter of Assumption \ref{noise}.
In  all  our schemes,  the error  can be splitted  in space  error  $err_0$ and the time error $err_1$.
The space error $err_0$ can be estimated as in Lemma \ref{serror} using results  from the proof of \cite[Theorem 4.1]{manto}.
 More precisely in the estimation of the poisson term in their case, we replace $\beta\in [0,2]$ by $\eta=2H+\beta-1$ defined in \assref{assumption5}. The time error  $err_1$ will be here splitted in three terms. 
 More precisely the deterministic $I_1$ related to the nonlinear function $F$, the fBm term $I_2$ and the Poisson term $I_3$.
 The estimation of  $I_1$ is done with the aid of Assumption \ref{der} and \eqref{mild_sol5} similarly  as the work done in \cite[(77)-(88)]{manto} for implicit and exponential schemes. As in the proof of \cite[Theorem 10]{Muka}, we use the Taylor expansion in Banach space (see \cite[(77)]{manto}) to estimate $I_1$ for SERS. The fBm term $I_2$ is done exactly as in the previous section for the  scheme without Poisson. By replacing $\beta\in [0,2]$ by $\eta=2H+\beta-1$ defined in \assref{assumption5}, the estimation of the  Poisson term $I_3$ is done as in \cite[Theorem 5.2]{manto} for implicit and SETD1 schemes using  Burkholder-Davis-Gundy Inequality, the  work  in \cite[Theorem 10]{Muka} and  preparatory results Lemma \ref{lemSERS1}-\ref{lemSERS4} for SERS.
\section{Numerical simulations}
\label{numerik}
\begin{figure}[!ht]
  \subfigure[]{
    \label{FIGIa}
    \includegraphics[width=0.47\textwidth]{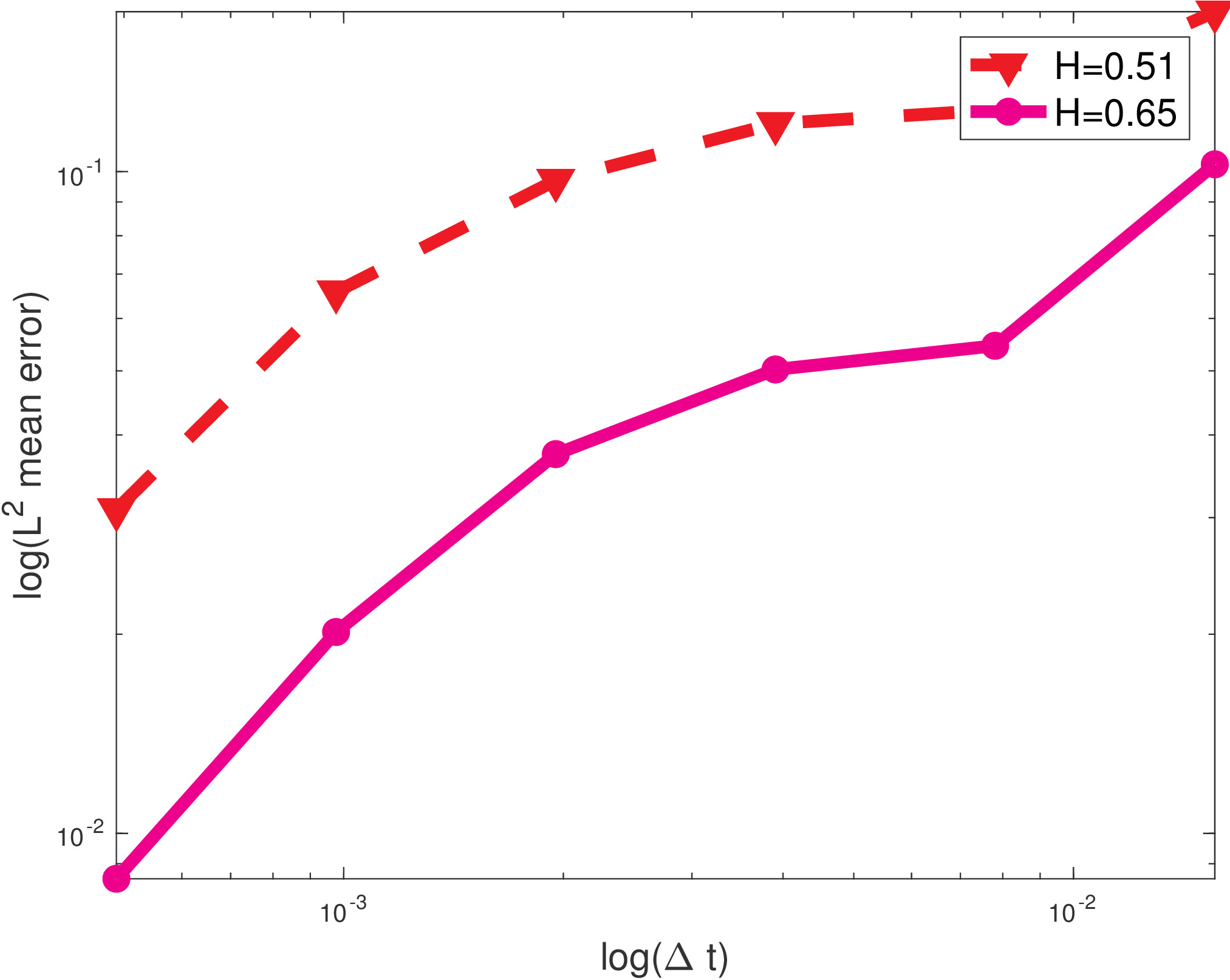}}
  \subfigure[]{
    \label{FIGIb}
    \includegraphics[width=0.47\textwidth]{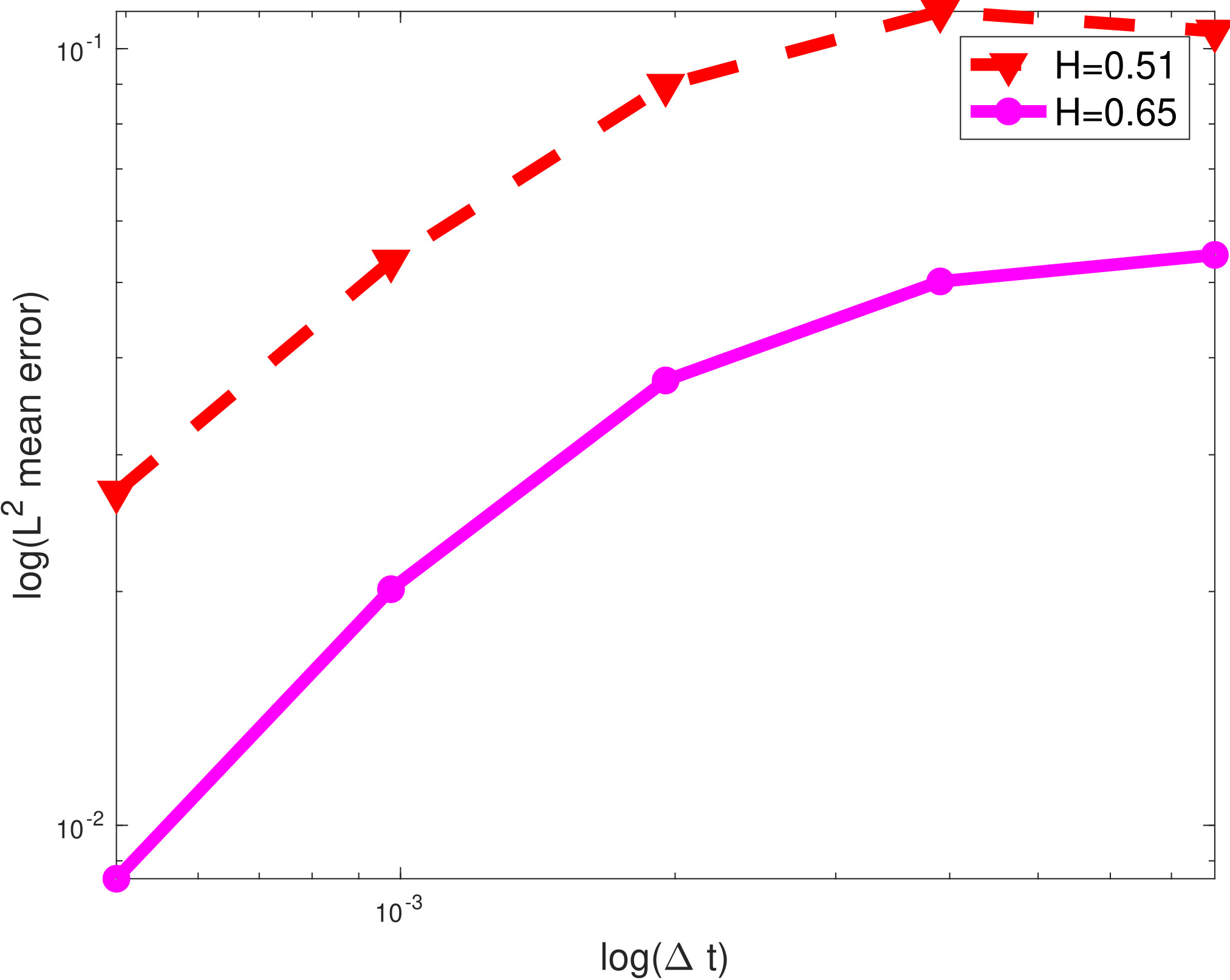}}
    \subfigure[]{
    \label{FIGIc}
    \includegraphics[width=0.47\textwidth]{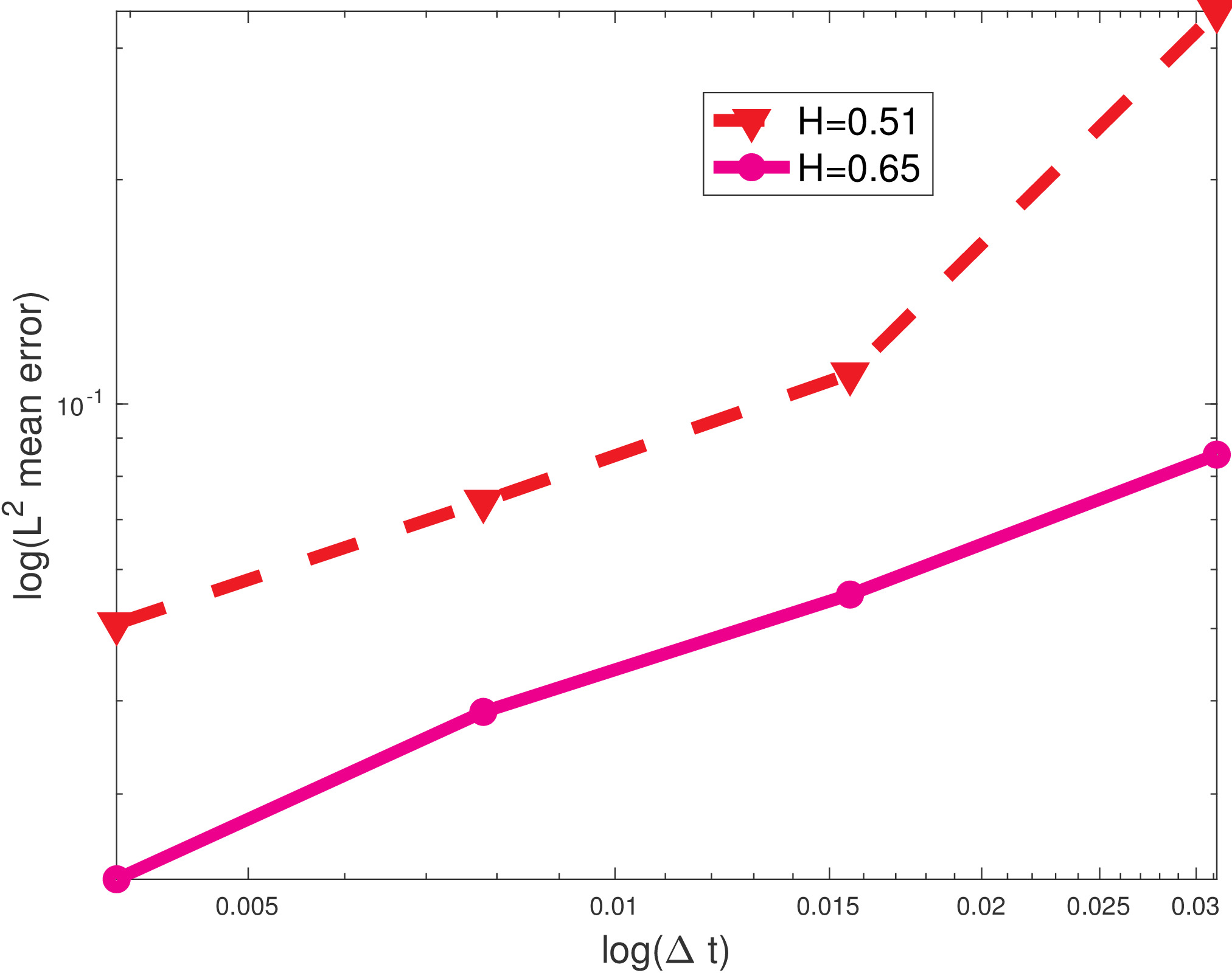}}
  \subfigure[]{
   \label{FIGId}
    \includegraphics[width=0.47\textwidth]{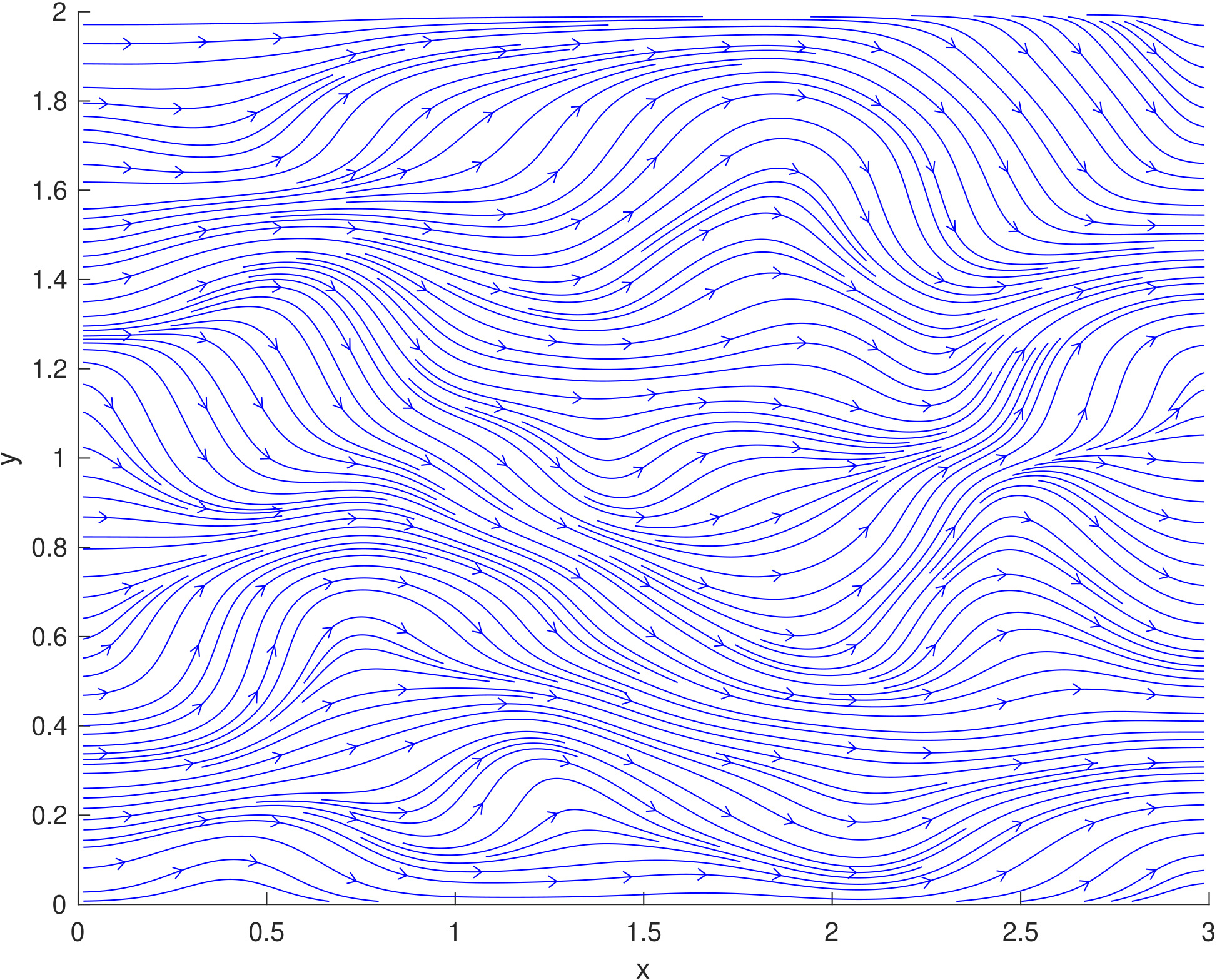}}
  \caption{Convergence in the root mean square $L^{2}$ norm at $T=1$ as a
    function of $\Delta t$ for implicit scheme (a), exponential  scheme (b) and exponential Rosenbrock scheme  (c). 
     We have  used here 50 realizations.  The streamline of the velocity field $\mathbf{q}$. } 
  \label{FIGI} 
  \end{figure}

In opposite  to  the standard   Brownian motion  where  the simulation is obvious,  the simulation of fBm  is not obvious and  is  an important research field in numerical analysis.  Keys methods for simulations of fBm are  Cholesky method \cite{Kam} and the circulant method \cite{G}, which will be used in this work to generate the fBm.  Here  we  consider the stochastic advection
diffusion reaction SPDE \eqref{pb1}-\eqref{secondorder} with  constant diagonal diffusion tensor  $ \mathbf{D}= 10^{-2} \mathbf{I}_2=(D_{i,j}) $ in \eqref{operator},
and mixed Neumann-Dirichlet boundary conditions on $\Lambda=[0,L_1]\times[0,L_2]$. 
The Dirichlet boundary condition is $X=1$ at $\Gamma=\{ (x,y) :\; x =0\}$ and 
we use the homogeneous Neumann boundary conditions elsewhere.
The eigenfunctions $ \{e_{i,j} \} =\{e_{i}^{(1)}\otimes e_{j}^{(2)}\}_{i,j\geq 0}
$ of  the covariance operator $Q$ are the same as for  Laplace operator $-\varDelta$  with homogeneous boundary condition  given by 
\begin{eqnarray*}
e_{0}^{(l)}(x)=\sqrt{\dfrac{1}{L_{l}}},\qquad 
e_{i}^{(l)}(x)=\sqrt{\dfrac{2}{L_{l}}}\cos\left(\dfrac{i \pi }{L_{l}} x\right),
\, i \in \mathbb{N}
\end{eqnarray*}
where $l \in \left\lbrace 1, 2 \right\rbrace,\, x\in \Lambda$. In the noise representation  \eqref{def_fBm}, we  have used
\begin{eqnarray}
\label{noise2}
 \lambda_{i,j}=\left( i^{2}+j^{2}\right)^{-(\beta +\delta)}, \, \beta>0,
\end{eqnarray} 
 for some small $\delta>0$.   We have  used  $b(x,t)=2$ in \eqref{secondorder}, so $\phi$ in \assref{noise} is obviously satisfied for $\beta=(0,1]$. 
 In our simulations, we have used $\delta=0.001$.
The function $f$ used in \eqref{nemystskii} to be  $f(x,z)= \frac{z}{1+z}$ for all $(x,z)\in \Lambda\times\mathbb{R}$. Therefore the corresponding Nemytskii operator   $F$ defined by \eqref{nemystskii} obviously satisfies 
\assref{non}.
We obtain the Darcy velocity field $\mathbf{q}=(q_i)$  by solving the following  system
\begin{equation}
  \label{couple1}
  \nabla \cdot\mathbf{q} =0, \qquad \mathbf{q}=-\mathbf{k} \nabla p,
\end{equation}
with  Dirichlet boundary conditions on 
$\Gamma_{D}^{1}=\left\lbrace 0,L_1 \right\rbrace \times \left[
  0,L_2\right] $ and Neumann boundary conditions on
$\Gamma_{N}^{1}=\left( 0,L_1\right)\times\left\lbrace 0,L_2\right\rbrace $ such that 
\begin{eqnarray*}
p&=&\left\lbrace \begin{array}{l}
1 \quad \text{in}\quad \left\lbrace 0 \right\rbrace \times\left[ 0,L_2\right]\\
0 \quad \text{in}\quad \left\lbrace L_1 \right\rbrace \times\left[ 0,L_2\right]
\end{array}\right. 
\end{eqnarray*}
and $- \mathbf{k} \,\nabla p (\mathbf{x},t)\,\cdot \mathbf{n} =0$ in  $\Gamma_{N}^{1}$. 
Note that $\mathbf{k}$ is the permeability tensor and $p$ the presure.  We use a random permeability field as in \cite[Figure 6]{Antofirst}.
 The streamline of the velocity field $\mathbf{q}$ are given in \figref{FIGI}(d). To deal with high  P\'{e}clet number,  we discretise in space using
finite volume method, viewed as a finite element method (see \cite{Antonio3}).
We take $L_1= 3$ and $L_2=2$ and  our reference solutions samples are numerical solutions using at time step of $\Delta t = 1/ 4096$. The errors are computed at the final time $T=1$.
The initial  solution  is $X_0=0$, so we can therefore expect high orders convergence, which depend  only on the noise term  and $H$.

\figref{FIGI}(a) is  the errors graph for the implicit scheme with different values of  $H$. We have observed
  that the order of convergence is $0.48$  in time for $ H=0.51$ and $\beta=1$,  $ 0.6476$ for $H=0.65$ and $\beta=1$ .
  
  \figref{FIGI}(b) is  the errors graph for the exponential scheme with two values of  $H$.  We have observed
  the order of convergence is $0.5012$ in time  for $ H=0.51$ and $\beta=1$,  $ 0.6653$ for $H=0.65$ and $\beta=1$.
  
\figref{FIGI}(c) is  the errors graph for the exponential Rosenbrock  scheme with two values of  $H$.  We have observed
  the order of convergence is $0.5562$ in time  for $ H=0.51$ and $\beta=1$,  $ 0.6197$ for $H=0.65$ and $\beta=1$. 

  As  we can observe,  our  numerical orders in time  are  close  to  our theoretical results in \thmref{cstrongconvthm} even if we have only used 50 samples in our Monte Carlo simulations.
  
\begin{figure}[!ht]
\begin{center}
  \subfigure[]{
    \label{FIGIIa}
    \includegraphics[width=0.47\textwidth]{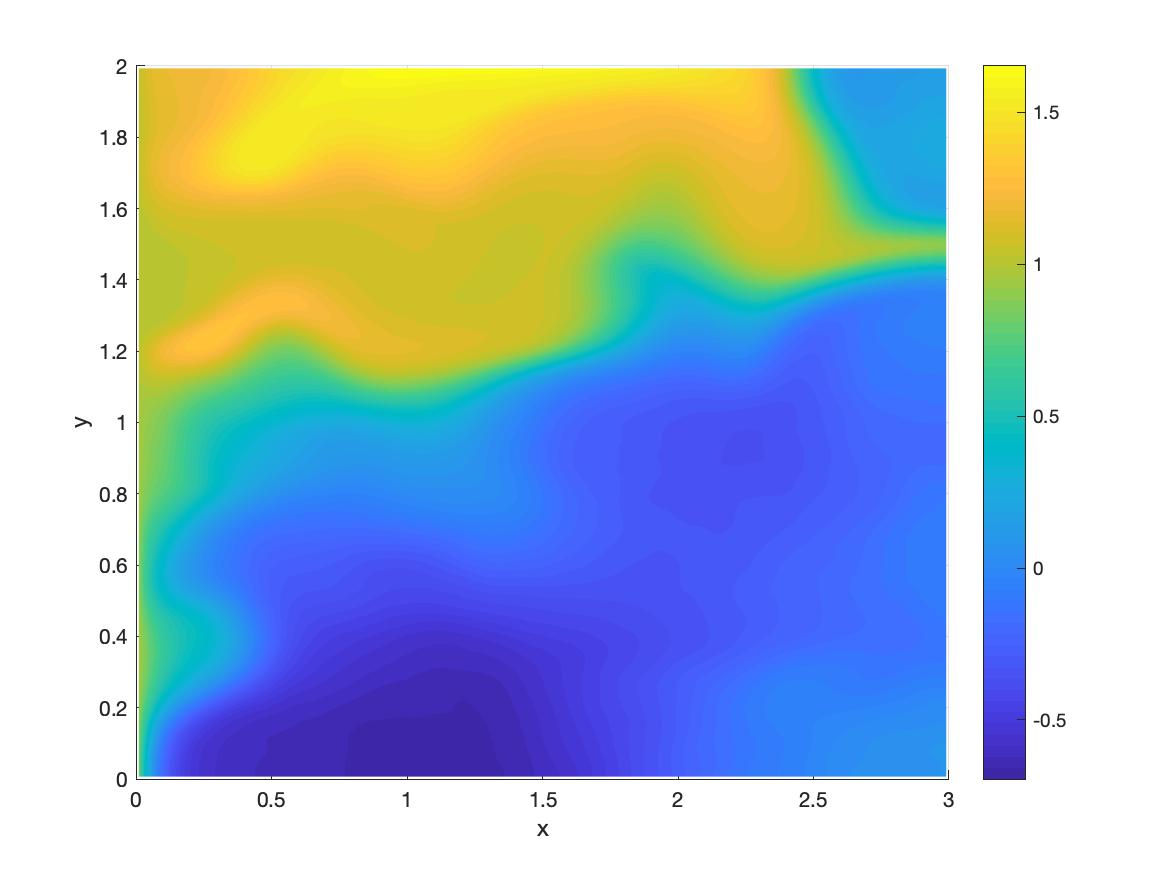}}
  \subfigure[]{
    \label{FIGIIb}
    \includegraphics[width=0.47\textwidth]{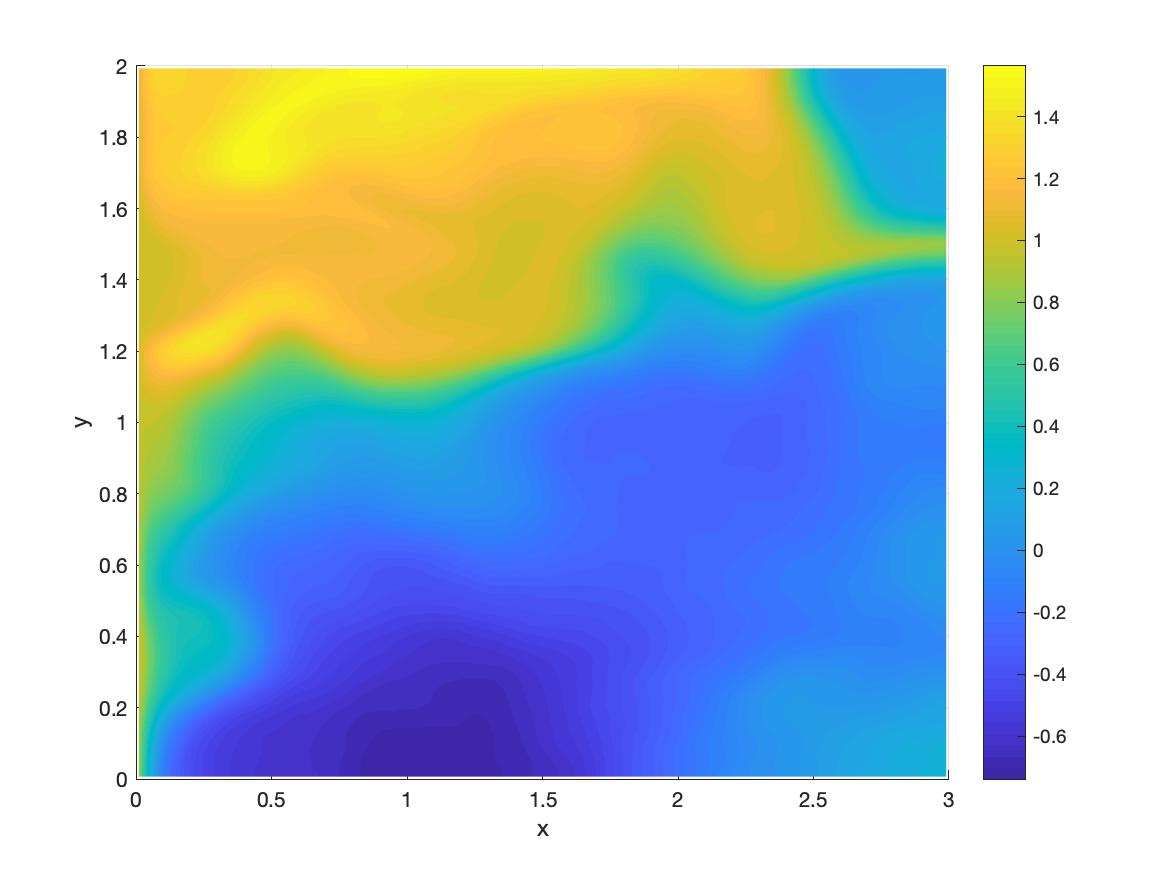}}
  \caption{Samples of the solution for $H=0.75$  (a) and $H=0.51$(b) } 
  \label{FIGII} 
  \end{center}
  \end{figure}

  \figref{FIGII} shows two samples of the solution for $H=0.75$ and $H=0.51$.  Here we have fixed $\beta=1$ and same  Gaussian  randoms numbers have been
   to generate our fBm samples. As we can observe, the   parameter $H$ has significant influence on the sample of the numerical solution.
  This is independent of our timestepping methods.
  
  \newpage
 
\section*{Acknowledgement}
 Aurelien Junior Noupelah 
thanks Prof  Louis   Aime  Fono  and  Prof  Jean Louis Woukeng for their  constant supports. 
We would like to thank  Jean Daniel Mukam  for very useful discussions.
\end{document}